%% file: padi9.tex
\documentstyle[eqsection]{article}

\include{macro}

\newcommand{\wG}{\widetilde{G}}

\newcommand{\dps}{\displaystyle}

\begin{document}

\title{{\bf Arnold's Diffusion in
nearly integrable isochronous Hamiltonian systems}}

\author{Massimiliano Berti and Philippe Bolle}
\date{}
\maketitle

{\bf Abstract:}
We consider the problem of Arnold's diffusion for
nearly integrable isochronous Hamiltonian systems.
We prove a shadowing theorem which improves the known
estimates for the diffusion time. We also develop
a new method for measuring the splitting of the
separatrices. As an application we justify, for three time scales systems,
that the splitting
is correctly predicted by the Poincar\'e-Melnikov function.
\footnote{Supported by M.U.R.S.T. Variational Methods and Nonlinear
Differential Equations.}
\\[2mm]
Keywords: Arnold's diffusion, shadowing theorem, splitting of separatrices,
heteroclinic orbits, variational methods.


\section{Introduction}

Through this paper we consider nearly integrable
isochronous Hamiltonian systems as
\be\label{eq:Hamg}
{\cal H}_\mu = \om \cdot I + \frac{p^2}{2} + (\cos q - 1) +
\mu f( \vphi, q ),
\ee
where $( \vphi, q ) \in {\bf T}^n \times {\bf T}^1 
:= ({\bf R}^n / 2\pi {\bf Z}^n) \times ({\bf R}/ 2\pi {\bf Z})$ are 
the angle variables,  $(I, p) \in {\bf R}^n \times {\bf R}^1$
are the action variables and $ \mu \geq 0 $ is a small real  
parameter. Hamiltonian $ {\cal H}_\mu $ describes a
system of $n$ isochronous harmonic oscillators
of frequencies $ \om $ weakly coupled with a pendulum.

When $ \mu = 0 $ the energy $ \om_i I_i $ of each oscillator is a
constant of the  motion.
The problem of {\it Arnold's diffusion} is whether,
for $ \mu \neq 0 $, there exist motions whose net effect is to transfer
energy from one oscillator to others.
This problem has been broadly investigated by many
authors also for non-isochronous systems, see for example \cite{Arn},
\cite{Bs}, \cite{BCV}, \cite{CG} and \cite{Xi}. In this paper we focus
on isochronous systems for which, in order to exclude trivial drifts
of the actions due to resonance phenomena,
it is standard to assume a diophantine condition
for the frequency vector $ \om $. Precisely we will always suppose
\begin{itemize}
\item
$(H1)$ $\exists \ \gamma > 0 $, $ \tau > n $
such that $| \om \cdot k| \geq \gamma / |k|^{\tau} $,
$ \forall k \in {\bf Z}^n, k \neq 0$.
\end{itemize}

The existence of Arnold's diffusion is usually proved
following the mechanism proposed in \cite{Arn}.
First one remarks that, for $ \mu = 0 $, Hamiltonian ${\cal H}_\mu$ admits
a continuous family of
$n$-dimensional partially hyperbolic invariant tori
${\cal T}_{I_0} = \{ (\vphi , I, q,p) \in  {\bf T}^n
\times  {\bf R}^n \times {\bf T}^1 \times {\bf R}^1 \ |
\  I = I_0, \ q = p = 0 \}$
possessing stable and unstable manifolds
$ W^s ({\cal T}_{I_0}) = W^u({\cal T}_{I_0}) =
\{ (\vphi , I, q,p) \in  {\bf T}^n
\times  {\bf R}^n \times {\bf T}^1 \times {\bf R}^1 \  |
  \ I = I_0, \ p^2 / 2 + ( \cos q - 1) = 0 \}.$
Then Arnold's mechanism is based on the following three main steps.
\begin{description}
\item {$Step \ (i)$} To prove that, for $ \mu $ small enough,
the perturbed stable and unstable manifolds $ W^s_\mu ({\cal T}_{I_0}^\mu )$
and $ W^u_\mu ({\cal T}_{I_0}^\mu) $ split and intersect transversally
(``splitting of the separatrices'');
\item {$Step \ (ii)$}
To prove the existence of a chain of tori connected by heteroclinic orbits
(``transition chain'');
\item {$Step \ (iii)$}
To prove, by a shadowing type argument, the existence
of an orbit such that the action variables
$I$ undergo a variation of $ O (1) $
in a certain time $ T_d $ called the {\it diffusion time}.
\end{description}
We point out that for isochronous systems assumption
$(H1)$ implies that {\it all} the unperturbed tori
${\cal T}_{I_0}$, with their
stable and unstable manifolds, persist, for $\mu$ small enough,
being just sligthly deformed. For this reason the construction of
the ``transition chain'' of step $(ii)$
is a straightforward consequence of step $(i)$.
This also happens, for the peculiar   choice of the perturbation,
in the non-isochronous system considered in \cite{Arn}.
  On the other hand, this is not the case
for general non-isochronous systems where the surviving perturbed
tori are separated by the gaps appearing in KAM constructions,
making the existence of chains of tori a difficult matter, see \cite{CG}.
We quote paper \cite{Xi} for a somewhat different mechanism of diffusion
where step $(ii)$ is bypassed using Mather's theory.
\\[2mm]
\noindent
In the present paper we address, 
for isochronous systems, the following two main questions
\begin{description}
\item{$1)$} Shadowing theorems and estimates of the diffusion time;
\item{$2)$} Splitting of separatrices.
\end{description}

Problem $1)$ has been intensively studied in the last years, see
for example \cite{BCV},\cite{CG},\cite{C1},\cite{G1} and \cite{M}
(we underline that \cite{BCV},\cite{CG},\cite{C1} and \cite{M}
deal also with non-isochronous systems).
Our general shadowing theorem
({\it thm.\ref{thm:main}-thm.\ref{thm:main1g}})
improves -for isochronous systems- the estimates
on the diffusion time obtained in the 
forementioned papers.

The estimate on the diffusion time
that we obtain (see expression (\ref{timediff})),
once   it  is verified that the stable and the unstable manifolds split,
is roughly the following:
the diffusion time $ T_d $ is estimated by the product
of the number of heteroclinic transitions $k$ ($=$ number of tori
forming the transition chain = heteroclinic jump/splitting)
and of the time $ T_s $ required for a single transition,
namely $T_d = k T_s$.
The time for a single transition $ T_s $  is bounded by the maximum
time between the   ``ergodization time'' of the torus 
${\bf T}^n $ run by the linear flow $ \om t $, 
and the time needed  to ``shadow''
homoclinic orbits for the quasi-periodically forced pendulum.

In order to highlight the improvement of our estimate of the
diffusion time let us consider 
the particular case of  ``a-priori unstable'' systems, 
i.e. when the frequency vector $ \om $ is considered as a constant 
independent of any parameter. 
In such a case
it is easy to evaluate, using the classical Poincar\'e-Melnikov theory,
that the splitting of the separatrices is $O(\mu)$. Then
our shadowing theorem yields the estimate for the diffusion time
$ T_d = O(  (1/ \mu ) \log (1 / \mu)) $, 
see {\it thm. \ref{thm:apu}-thm.\ref{apug}}.

Such estimate answers to a question raised in
\cite{Lo} (sec.7) proving that, at least for isochronous systems,
it is possible to reach the maximal speed of diffusion
$ \mu /| \log \mu |$.
On the contrary the estimate on the diffusion time obtained
in \cite{CG} is $ T_d >> O( \exp{(1/\mu)} )$
and is improved in \cite{G1} to be $ T_d = O(\exp{(1/ \mu )}) $.
Recently in \cite{BCV} by means of Mather's theory
the estimate on the diffusion time has been improved to be
$ T_d = O ( 1/ \mu^{2 \t +1} ) $.
In \cite{C1} it is obtained via geometric methods
that $ T_d = O ( 1/ \mu^{\t +1} ) $.
It is worth pointing out that the estimates given in  
\cite{BCV} and \cite{C1}, which yet provide a diffusion time
polinomial in the splitting,
depend on the diophantine exponent $ \tau $ and hence on 
the number of rotators $ n $.
On the contrary our estimate is independent of $ n $.

The main reason for which we are able to improve also the estimates
of \cite{BCV} and \cite{C1} is that our shadowing orbit
can be chosen, at each transition, to approach the
homoclinic point, only up to a distance $O(1)$
and not $O(\mu)$ like in \cite{BCV} and \cite{C1}.
This implies that the time spent by our diffusion orbit
at each transition is $T_s = O( \log (1/\mu) )$. Since
the number of tori forming the transition chain is equal to
$ O(1/ {\rm splitting}) = O(1/ \mu)$ the diffusion time is finally
estimated by $ T_d = O( (1/ \mu) \log (1 / \mu)) $.
\\[1mm]
\indent
Regarding the method of proof, we use a variational technique 
inspired by \cite{BB} and \cite{BBS}. One advantage of this approach 
is that the same arguments can be also used when the hyperbolic part is
a general Hamiltonian in ${\bf R}^{2m}$, $ m \geq 1 $, possessing
one hyperbolic equilibrium and a transversal homoclinic orbit.
Nevertheless we have developed all the details in 
the case that the hyperbolic part is
the standard one-dimensional pendulum because it is
the model equation to study Arnold's diffusion near a simple-resonance.

Furthermore we also remark that our proof of 
theorem \ref{thm:main1g} is completely self-contained
in the sense that, unlike the known approaches 
(excepted \cite{Xi}), we do not make use of any KAM-type result
for proving, under assumption $(H1)$, the persistence of
invariant tori, see {\it thm. \ref{thm:invto}}.
\\[2mm]
\indent
In sections 4-5 we study problem 2).
Detecting the splitting of the separatrices
becomes a very difficult problem when the frequency vector $ \om = \om_\e $
depends on some small parameter $ \e $ and
contains some ``fast frequencies'' $ \om_i = O( 1/ \e^b )$, $ b > 0 $.
Indeed, in this case, the oscillations of
the Melnikov function along some directions turn out to be exponentially
small with respect to $ \e $ and then
the naive Poincar\'e-Melnikov expansion
provides a valid measure of the splitting only for
$ \mu $ exponentially small with respect to $ \e $.
Much literature in the last years has been
devoted to overcome this problem,
see for example \cite{DGTJ},\cite{GGM}, \cite{GGM1}, \cite{LMS}
and \cite{Sa}.
In the present paper, in order to justify the dominance
of the Poincar\'e-Melnikov function when $ \mu = O(\e^p )$,
we extend the approach originally used in \cite{An1} for dealing with
rapidly periodic forced systems.
Up to a change of variables close to the identity, we prove
(see {\it thm. \ref{thm:tfc}}) an exponentially
small upper bound for the Fourier coefficients of the splitting.
As an application we provide some results
({\it thm. \ref{thm:tts}, thm. \ref{lem:3ts})} on
the splitting of the separatrices and the diffusion time
({\it thm. \ref{thm:3tsm}}) for three time scales systems
$$
{\cal H}_\e =  \frac{1}{\sqrt{\e}} I_1 + \e^a \b \cdot I_2 +
\frac{p^2}{2} + ( \cos q - 1)  +
\mu  ( \cos q - 1) f( \vphi ),  \
I_1 \in {\bf R}^1, I_2, \b \in {\bf R}^{n-1}, \ n \geq 2,
$$
for $ \mu \e^{-3/2} $ small.
This improves the main theorem I in \cite{PV} which holds for
$ \mu = \e^p $, $ p > 2 + a $.
With respect to \cite{GGM}, which deals for more general systems,
we remark that our results hold in any dimension,
while the results of \cite{GGM} are proved for 2 rotators only.

Theorem \ref{thm:tfc} is also the starting point
for proving the splitting of the separatrices in presence of
two high frequencies, assuming as in
\cite{DGTJ},\cite{GGM1}, \cite{LMS} and \cite{Sa} suitable
 hypotheses on the perturbation term. We do not address
this problem in this paper.

On the other hand theorem \ref{thm:tts} is the starting point to prove,
for $ n \geq 3 $, the existence of diffusion solutions such that
the action variables $ I_2 $ undergo a variation $ O(1) $ in polynomial time
(while the $ I_1 $ action variable does not change considerably).
This phenomenon can not be deduced by the estimates,
given in \cite{GGM} and \cite{PV}, on the ``determinant of the splitting''
which does not distinguish among ``slow'' and ``fast'' directions
and would give rise to exponentially large diffusion times.
This type of results are contained in the forthcoming paper \cite{BB3}.
\\[2mm]
\indent
The paper is organized as follows:
in section 2 we prove the   shadowing theorem when the perturbation term is
$ f (\vphi, q) = ( 1- \cos q) f ( \vphi ) $.
In section 3 we show how to prove the theorem for general
perturbation terms $f ( \vphi, q )$. In section 4 we provide
the theorem on the Fourier coefficients of the splitting
and in section 5 we consider three   time scales systems.
\\[1mm]
The results of this paper have been announced in \cite{BBN}.
\\[2mm]
After this paper was completed we learned by prof. Bolotin 
about the recent preprint \cite{Tr} which deals with a-priori 
unstable Hamiltonian systems time periodically forced.
Among many results, in theorem 2 of \cite{Tr} a shadowing theorem 
for a symplectic separatrix map is proved 
providing an estimate on the diffusion speed
$ \mu /| \log \mu |$, like our. However theorem 2 applies just to an
approximation of the symplectic map describing the true dynamics
of the system and, moreover, requires the hyperbolic part to 
be just two dimensional.
\\[2mm]
{\bf Acknowledgments:} The first author wishes to thank Prof. G. Gallavotti
for stimulating discussions.

\section{The shadowing theorem}

We first develop our approach when the perturbation term
$ f( \vphi, q )= (1 - \cos q ) f( \vphi )$ so that the tori
$ {\cal T}_{I_0} $  are still invariant for $ \mu \neq 0 $.
The equations of motion derived by Hamiltonian ${\cal H}_\mu$ are
\be\label{eqmotion}
\dot{\vphi} = \om, \qquad
\dot{I} = - \mu (1 - \cos q ) \ \partial_{\vphi}f(\vphi ), \qquad
\dot{q} = p, \qquad
\dot{p} = \sin{q} - \mu \sin q \ f (\vphi ).
\ee
The dynamics on the angles $ \vphi $ is given by
$ \vphi (t) = \om t + A $   so that (\ref{eqmotion})
  are reduced  to the quasi-periodically forced pendulum
\be\label{pendper}
- \ddot{q} + \sin{q} = \mu \ \sin q \ f ( \om t + A)
\ee
corresponding to the Lagrangian
\be\label{lagraper}
{\cal L}_\mu (q, \dot{q},t) = \frac{{\dot q}^2}{2} + (1- \cos q) +
\mu (\cos q - 1)f(\om t + A).
\ee
For each solution $ q(t) $ of (\ref{pendper}) one recovers the dynamics
of the actions $ I (t) $ by quadratures in (\ref{eqmotion}).

\subsection{$1$-bump homoclinic and heteroclinic solutions}

For $ \mu = 0 $ equation (\ref{pendper})
is autonomous and possesses the homoclinic (mod. $2\pi$) solutions
$ q_\teta (t) = 4 \ {\rm arctg} (\exp {(t - \teta)})$,
$ \teta \in {\bf R}$.
Using the Contraction Mapping Theorem we now prove that, near
the unperturbed homoclinic solutions $ q_\teta (t)$, there exist,
for $ \mu $ small enough, {\it ``pseudo-homoclinic solutions''}
$ q_{A,\teta}^\mu (t)$ of equation (\ref{pendper}).
$ q_{A,\teta}^\mu (t)$ are true solutions of
(\ref{pendper}) in $( - \infty, \teta)$ and $(\teta, +\infty)$;
at time $ t = \teta $ such pseudo-solutions are glued with continuity
at value $ q_{A,\teta}^\mu ( \teta ) = \pi  $
and for $ t \to \pm \infty $ are asymptotic to the equilibrium $0$ mod $2
\pi$.
\begin{lemma}\label{lem:1bump}
There exist $ \mu_0, C_0 > 0 $ such that $\forall 0 < \mu  \leq \mu_0 $,
$ \forall \om \in {\bf R}^n $, $ \forall \teta \in {\bf R} $,
there exists a unique function
$ q_{A, \teta}^\mu (t) : {\bf R} \to {\bf R}$, smooth
in $ ( A, \teta, \mu )$, such that
\begin{itemize}
\item{$(i)$}
$ q_{A, \teta}^\mu (t)$ is a solution of (\ref{pendper})
in each interval $(- \infty, \teta)$ and $(\teta, +\infty)$
and $ q^\mu_{A, \teta} ( \teta ) = \pi $;
\item{$(ii)$}
$\max \Big( | q^\mu_{A, \teta }(t) - q_\teta (t)|,
|{\dot q}^\mu_{A, \teta }(t) - {\dot q}_\teta (t)| \Big) \leq C_0 \mu 
\exp({- \frac{ |t - \teta |}{2}})$, $ \forall t \in {\bf R}$;
\item{$(iii)$}
$ q^\mu_{A, \teta } (t) = q^\mu_{A + k 2 \pi, \teta } (t), \
\forall k \in {\bf Z}^n $;
\item{$(iv)$}
$q_{A, \teta + \eta }^\mu (t + \eta ) = q_{A + \om \eta, \teta}^\mu (t), \
\forall \teta, \eta \in {\bf R}$.
\item{$(v)$} 
$ \max \Big( | \partial_A  q^\mu_{A, \teta } (t) |,
| \partial_A \dot{q}^\mu_{A, \teta } (t)|,
|\om \cdot \partial_A  q^{\mu}_{A, \teta } (t) |,
|\om \cdot \partial_A \dot{q}^{\mu}_{A, \teta } (t)| \Big) =  
O \Big(\mu  \exp(- \frac{ |t - \teta |}{2}) \Big) $. 
\end{itemize}
\end{lemma}

\begin{pf}
Proof in the appendix. Note that in $(v)$ the bound 
of $|\om \cdot \partial_A  q^\mu_{A, \teta } (t) |,
|\om \cdot \partial_A \dot{q}^\mu_{A, \teta } (t)|$ is uniform in
$\om$.
\end{pf}

We can then define the function
$ F_\mu : {\bf T}^n \times {\bf R} \to {\bf R} $
as the action functional of Lagrangian (\ref{lagraper})
evaluated on the ``1-bump pseudo-homoclinic solutions''
$ q_{A, \teta }^\mu (t) $, namely
\be\label{def:homo}
F_\mu (A, \teta) =
\int_{ - \infty}^\teta {\cal L}_\mu ({  q}_{A,\teta}^\mu (t),
\dot{  q}^\mu_{A,\teta} (t), t) \ dt
+ \int_\teta^{+ \infty}
{\cal L}_\mu ({  q}_{A,\teta}^\mu (t),
\dot{  q}^\mu_{A, \teta } (t),t) \ dt
\ee
and the {\it ``homoclinic function''} $ G_\mu: {\bf T}^n \to {\bf R} $ as
\be\label{eq:homospli}
G_\mu ( A ) = F_\mu ( A, 0 ).
\ee
Since $ q_{A,\teta}^\mu (t) $ converges
exponentially fast to $0$, mod $2 \pi$, the integrals
in (\ref{def:homo}) are convergent.
Note that the homoclinic function $G_\mu $ is independent of $I_0$.
By property $(v)$ of lemma \ref{lem:1bump}
the following invariance property holds
$$
F_\mu (A, \teta + \eta) = F_\mu (A + \om \eta, \teta), \quad \forall
\teta, \eta \in {\bf R},
$$
and in particular
\be\label{eq:invaspe}
F_\mu  (A , \teta ) = G_\mu (A + \om \teta), \quad \forall
\teta \in {\bf R}.
\ee

\begin{remark}
The homoclinic function $ G_\mu $ is the difference
between the generating functions
$ {\cal S}_{\mu,I_0}^\pm (A, q_0 )$ of the stable and the unstable manifolds
$ W_\mu^{s,u} ( {\cal T}_{I_0} )$
(which in this case are {\it exact} Lagrangian manifolds)
at section $ q_0 = \pi $, namely
$G_\mu (A) = {\cal S}_{\mu, I_0}^- (A, \pi ) -
{\cal S}_{\mu, I_0}^+ ( A, \pi ) $.
Indeed can be easily verified that
$$
{\cal S}_{\mu,I_0}^+ (A,q_0) = I_0 \cdot A - \int_0^{+ \infty}
\frac{({ \dot q }^\mu_{A, q_0} (t) )^2}{2} + ( 1 - \cos q^\mu_{A, q_0}(t)) +
\mu (\cos q^\mu_{A, q_0}(t) - 1) f( \om t + A ) \ dt,
$$
where $ q^\mu_{A, q_0} (t) $ is the unique solution
 of (\ref{pendper}) near $q_0 (t)$
with $ q^\mu_{A, q_0} (0) = q_0 $ and
$ \lim_{t \to + \infty} q^\mu_{A, q_0} (t) = 2 \pi $.
Analogously
$$ {\cal S}_{\mu, I_0}^- ( A, q_0) := I_0 \cdot A + \int_{- \infty}^0
\frac{({ \dot q }^\mu_{A, q_0} (t))^2}{2} + ( 1 - \cos q^\mu_{A, q_0}(t) ) +
\mu (\cos q^\mu_{A, q_0} (t)- 1) f( \om t + A ) \ dt,
$$
where $ q^\mu_{A, q_0} (t) $ is the unique solution
of (\ref{pendper}) near $ q_0 (t) $ with $ q^\mu_{A, q_0} (0) = q_0 $ and
$ \lim_{t \to - \infty} q^\mu_{A, q_0} (t) = 0 $.
\end{remark}

\begin{lemma}\label{criteta}
The derivative of $ \teta \to F_\mu (A, \teta) $ satisfies
\be
\partial_\teta F_\mu (A, \teta) =
\frac{({{\dot q}^\mu_{A,\teta}})^2 ( \teta^+ )}{2} -
\frac{({{\dot q}^\mu_{A,\teta}})^2 ( \teta^- )}{2}.
\ee
\end{lemma}

\begin{pf}
There holds
\begin{eqnarray*}
\partial_\teta F_\mu (A, \teta) & = &
\frac{({{\dot q}^\mu_{A,\teta}})^2 ( \teta^- )}{2} -
\frac{({{\dot q}^\mu_{A,\teta}})^2 ( \teta^+ )}{2} +
\int_{- \infty }^\teta \ldots \\
& + & \int_{ \teta }^{+ \infty}
{{\dot q}^\mu_{A,\teta}}( t) \partial_\teta {{\dot q}^\mu_{A,\teta}}(t)
+ \Big( \sin q_{A,\teta}^\mu (t) - \mu \sin q_{A,\teta}^\mu (t)
f ( \om t + A) \Big) \partial_\teta q_{A,\teta}^\mu (t) \ dt.
\end{eqnarray*}
Integrating by parts and
using also that $ q_{A,\teta}^\mu (t) $ solves (\ref{pendper}), we obtain
\be\label{eq:pezzetti}
\partial_\teta F_\mu ( A, \teta ) =
\frac{1}{2} ({{\dot q}^\mu_{A,\teta}})^2 ( \teta^- ) -
\frac{1}{2} ({{\dot q}^\mu_{A,\teta}})^2 ( \teta^+ ) +
\Big[ \partial_\teta q_{A,\teta}^\mu (t) {{\dot q}^\mu_{A,\teta}}( t)
\Big]_{-\infty}^{\teta^-} +
\Big[ \partial_\teta q_{A,\teta}^\mu (t) {{\dot q}^\mu_{A,\teta}}( t)
\Big]_{\teta^+}^{+\infty}.
\ee
Since $ \forall \teta \in {\bf R} $ $ q_{A,\teta}^\mu (\teta) = \pi $,
deriving in $ \teta $ we obtain
$ \partial_\teta q_{A,\teta}^\mu (\teta) +
{{\dot q}^\mu_{A,\teta}}( \teta) = 0 $; hence from (\ref{eq:pezzetti})
and using that
$ \lim_{t \to \pm \infty} {{\dot q}^\mu_{A,\teta}}(t) = 0 $ we deduce
lemma \ref{criteta}.
\end{pf}

By lemma \ref{criteta} if $ \partial_\teta F_\mu (A , \teta) = 0 $ then
$ q_{A, \teta}^\mu (t)$ is a true homoclinic (mod. $2\pi$)
solution of (\ref{pendper}).
Then, for each $ I_0 \in {\bf R}^n $,
\be\label{orbit}
\Big( \om t + A, I_\mu (t), q_{A, \teta }^\mu (t),
{\dot q}_{A, \teta}^\mu (t) \Big)
\ee
where
\be\label{salto}
I_\mu (t) = I_0 -  \mu \int_{- \infty}^t
(1 - \cos q^\mu_{A, \teta}( s ))
\partial_{\vphi} f (\om s + A) ds
\ee
is a solution of $ {\cal H}_\mu $ emanating at $ t = -\infty $ from
torus $ {\cal T}_{I_0} $.
Since $ q^\mu_{A, \teta} $ converges exponentially fast to
the equilibrium, the ``jump'' in the action variables
$ I_\mu ( + \infty ) - I_0 $ is finite.
We shall speak of homoclinic orbit to the torus
${\cal T}_{I_0}$ when the jump is zero, and of
heteroclinic from ${\cal T}_{I_0}$ to
${\cal T}_{I_{\mu}(+\infty)}$ when the jump is not zero.
Moreover the next lemma
says that such jump is given by $ \partial_A F_\mu (A, \teta ) $:

\begin{lemma}\label{lem:jump}
Let $ \partial_\teta F_\mu ( A, \teta ) = 0 $ then
$ I_\mu (t)$ given in (\ref{salto}) satisfies
\be\label{eq:jump}
\partial_A F_\mu (A, \teta ) =
\int_{ - \infty }^{+ \infty } {\dot I}_\mu (t) \ dt =
I_\mu ( + \infty ) - I_0  < + \infty .
\ee
In particular if $( A, \teta) $ is a critical
point of $F_\mu (A, \teta)$ then
(\ref{orbit}) in a homoclinic orbit to torus ${\cal T}_{I_0}$.
\end{lemma}

\begin{pf}
There holds
\begin{eqnarray*}
\partial_A F_\mu (A, \teta) & = &
\int_{- \infty}^{+\infty}
{{\dot q}^\mu_{A,\teta}}( t) \partial_A {{\dot q}^\mu_{A,\teta}}(t)
+ \sin q_{A,\teta}^\mu (t) \partial_A q_{A,\teta}^\mu (t) \\
& - & \mu \sin q_{A,\teta}^\mu (t) f(\om t + A) \partial_A q_{A,\teta}^\mu
(t)
- \mu ( 1 - \cos q_{A,\teta}^\mu (t)) \partial_\vphi f(\om t + A) \ dt.
\end{eqnarray*}
Integrating by parts, since $q_{A,\teta}^\mu (t)$ solves
(\ref{pendper}), and using that
$ \lim_{t \to \pm \infty} {{\dot q}^\mu_{A,\teta}}(t) = 0 $,
we deduce
\be\label{eq:inte}
\partial_A F_\mu (A, \teta) =
\int_{-\infty}^{+\infty}
- \mu ( 1 - \cos q_{A,\teta}^\mu (t) ) \partial_\vphi f( \om t + A ) \ dt.
\ee
We deduce from (\ref{eq:inte}) equality (\ref{eq:jump}).
\end{pf}

By the invariance property (\ref{eq:invaspe}) if $ B $
is a critical point of the homoclinic function $ G_\mu $, then,
 for all $(A, \teta )$ such that $ A + \om \teta = B$,
(\ref{orbit}) are homoclinic solutions to each torus
${\cal T}_{I_0}$. These homoclinics are not geometrically distinct since,
by the autonomy of ${\cal H}_\mu$, they are all obtained by time translation
of the same homoclinic orbit.
By the  Lusternik-Schirelman
category theory, since cat ${\bf T}^n = n+1 $, the function
$ G_\mu : {\bf T}^n \to {\bf R} $ has at least $ n+1 $ distinct
critical points. This  proves (see also \cite{LMS})

\begin{theorem}
Let $ 0 < \mu \leq \mu_0$. $ \forall I_0 \in {\bf R}^n $
there exist at least $ n+1 $ homoclinic orbits geometrically dinstict to
${ \cal T}_{I_0}$.
\end{theorem}

From the conservation of energy   a
heteroclinic orbit between ${\cal T}_{I_0}$ and ${\cal T}_{I_0'}$,
if any, must satisfy the energy relation
\be \label{eq:equalen}
\om \cdot I_0  = \om \cdot I_0'.
\ee
By lemma \ref{lem:jump}
a critical point of $ F_{\mu, I_0, I_0'}(A, \teta) $, defined by
$ F_{\mu, I_0, I_0'} (A, \teta)=$ $F_\mu (A, \teta) -
(I_0' - I_0) \cdot A =$ $G_\mu (A + \om \teta ) - (I_0' - I_0) \cdot A$,
gives rise to a heteroclinic solution joining the tori
$ {\cal T}_{I_0} $ to $ { \cal T}_{I_0'} $.
If the energy condition (\ref{eq:equalen}) holds then
the function $ F_{\mu, I_0, I_0'} (A, \teta) $ satisfies the
invariance property
\be
F_{\mu, I_0, I_0'} (A, \teta) = G_\mu ( A + \om \teta )
- (I_0' - I_0) \cdot (A + \om \teta)= G_{\mu, I_0, I_0'}
( A + \om \teta ).
\ee
where
\be\label{hetfunc}
G_{\mu, I_0, I_0'} ( B ) :=  G_\mu ( B ) - (I_0' - I_0) \cdot B.
\ee
Note that   $ G_{\mu, I_0, I_0'} $
is not $2\pi {\bf Z}^n$-periodic, and it
  might possess no critical point even for $ |I_0' - I_0| $ small.
However near   a  homoclinic orbit to ${\cal T}_{I_0}$ satisfying
some ``transversality condition''
there exist heteroclinic solutions connecting nearby tori
${\cal T}_{I_0'}$.   As an example,  the following
theorem holds, where $B_{\rho} (A_0)$ denotes an open ball in
${\bf R}^n$ (covering space of ${\bf T}^n$).

\begin{theorem}
Assume that there exist $ A_0 \in {\bf T}^n $, $\delta > 0 $ and
$ \rho > 0 $ such that
$\inf_{\partial B_\rho ( A_0 )} G_\mu > \inf_{B_\rho ( A_0 )} G_\mu +
\delta $. Then for all $ I_0, I_0' \in {\bf R}^n $ satisfying
$ ( I_0 - I_0') \cdot \om = 0 $ and $|I_0 - I_0'|
\leq \delta/ (2 \rho) $
there exists  a heteroclinic solution
of $ {\cal H}_\mu $ connecting ${\cal T}_{I_0}$ to ${\cal T}_{I_0'}$.
\end{theorem}

\subsection{The $k$-bump pseudo-homoclinic solutions}

We prove in the next lemma
the existence of pseudo-homoclinic solutions $ q_{A,\teta}^L (t)$
of the quasi-periodically forced pendulum (\ref{pendper}) which
turn $ k $ times along the separatrices
and are asymptotic to the equilibrium for $ t \to \pm \infty $.
Such pseudo-homoclinics $ q_{A,\teta}^L (t) $
are found, via the Contraction Mapping Theorem,
as small perturbations of a chain of ``1-bump pseudo-homoclinic
solutions'' obtained in lemma \ref{lem:1bump}.

\begin{lemma}\label{lem:kheter}
There  exist $ C_1, L_1 > 0 $ such that $\forall \om
\in {\bf R}^n $, $ \forall 0 <  \mu  \leq \mu_0 $,
$ \forall k \in {\bf N}$,
 $\forall  L > L_1  $, $ \forall \teta_1 < \ldots < \teta_k $ with
$ \min_i | \teta_{i+1} - \teta_i | > L$,
there exists a unique pseudo-homoclinic solution
$ q_{ A, \teta }^L ( t ): {\bf R} \to {\bf R} $, smooth
in $ ( A, \teta, \mu ) $ which is a true solution of (\ref{pendper})
in each interval $ (- \infty, \teta_1) $, $ ( \teta_i, \teta_{i+1} ) $
($ i = 1, \ldots, k-1 $), $ ( \teta_k, +\infty ) $ and
\begin{itemize}
\item{$(i)$}
$ q_{A, \teta }^L ( \teta_i ) = \pi ( 2 i - 1 ) $,
$ q_{A, \teta }^L (t) = q^\mu_{ A, \teta_1 } ( t ) $ in
$( - \infty, \teta_1 ) $ and
$ q_{A, \teta }^L ( t ) = 2 \pi (k-1) + q^\mu_{A, \teta_k }(t) $ in
$( \teta_k, + \infty )$;
\item{$(ii)$}
$|| q_{A, \teta }^L - q^\mu_{A, \teta_i}||_{W^{1, \infty }(J_i)} \leq
C_1 \exp({-C_1  L }) $
\ {\rm where} \ $ J_i = (\teta_i, (\teta_i + \teta_{i+1})/2 ),
\ \forall \ i = 1, \ldots, k-1$;
\item{$(iii)$}
$ || q_{A, \teta }^L  -
q^\mu_{A, \teta_{i+1} }||_{W^{1, \infty} (J_i')} \leq
C_1 \exp({- C_1 L}) $
\ {\rm where}\ $J_i' = ((\teta_i + \teta_{i+1})/2, \teta_{i+1}),
\ \forall \ i = 1, \ldots, k-1$;
\item{$(iv)$}
$ q^L_{A, \teta } (t) = q^L_{A + k 2 \pi, \teta } (t), \
\forall k \in {\bf Z}^n $;
\item{$(v)$}
$ q_{A, \teta + \eta }^L (t + \eta ) = q_{A + \om \eta, \teta}^L (t)$, \
$\forall \teta, \eta \in {\bf R}$.
\end{itemize}
\end{lemma}

\begin{pf}
In the appendix.
\end{pf}

We consider the Lagrangian action functional evaluated
on the pseudo-homoclinic solutions  $ q_{ A, \teta }^L $ given
by lemma \ref{lem:kheter}
depending on $ n + k $ variables
$$
F_\mu^k (A_1, \ldots, A_n, \teta_1, \ldots, \teta_k ) =
\int_{- \infty}^{\teta_1} {\cal L}_\mu (q_{A,\teta}^L (t),
{\dot{  q}}_{A,\teta}^L (t),t) \ dt \ +
$$
$$
\sum_{i=1}^{k-1} \int_{\teta_i}^{\teta_{i+1}}
{\cal L}_\mu (q_{A,\teta}^L (t),
{\dot q}_{A,\teta}^L (t),t) \ dt
+ \int_{\teta_k}^{+ \infty}
{\cal L}_\mu (q_{A,\teta}^L (t),
{\dot q}_{A, \teta}^L (t),t) \ dt.
$$
By lemma \ref{lem:kheter}-$v$ the following invariance
property holds
\be \label{invariance}
F_\mu^k ( A,\theta + \eta ) = F_\mu^k ( A + \eta \om, \teta ),
\qquad  \forall \teta, \eta \in {\bf R}.
\ee
Let $ {\cal F}_\mu^k : {\bf T}^n \times {\bf R}^k \to {\bf R} $
be   the ``$ k $-bump heteroclinic function'' defined by
\be\label{eq:kbhf}
{\cal F}_\mu^k ( A , \teta ) := F_\mu^k (A, \teta) - (I_0'- I_0) \cdot A.
\ee
Arguing as in lemma \ref{lem:jump} we have
\begin{lemma}\label{lem:heter}
$\forall I_0, I_0' \in {\bf R}^n $,
if $ ( A, \teta ) $ is a critical point of ${\cal F}_\mu^k (A, \teta)$, then
$ (\om t + A, I_\mu (t),$ $ q_{A, \teta}^L (t),
{\dot q}_{A, \teta }^L (t) )$
where
$ I_\mu (t) = I_0 - \mu \int_{-\infty}^t (1 -
\cos q_{ A, \teta }^L ( s ) ) \partial_\vphi f (\om s+ A)ds $
is   a heteroclinic solution connecting  $ {\cal T}_{I_0} $ to
$ {\cal T}_{I_0'} $.
\end{lemma}

By lemma \ref{lem:heter} we need to find
critical points of $ {\cal F}_\mu^k (A, \teta) $.
When $\min_i (\teta_{i+1} - \teta_i) \to + \infty $
the ``$k$-bump homoclinic function'' $ F_\mu^k ( A, \teta) $ turns out to be
well approximated simply by the sum of 
$ F_\mu ( A, \teta_i) $ according to the following lemma.
We set $ \teta_0 = - \infty$ and $ \teta_{k+1} = + \infty $.

\begin{lemma}\label{approxsum}
There   exist $ C_2, L_2 > 0 $ such that
$ \forall \om \in {\bf R}^n $, $ \forall 0<  \mu \leq \mu_0 $,
$ \forall L > L_2 $, $ \forall \teta_1 < \ldots < \teta_k $ with
$ \min_i ( \teta_{i+1} - \teta_i ) > L $
\be\label{eq:somma}
F^k_\mu (A, \theta_1, \cdots,
\theta_k ) = \sum_{i=1}^k F_\mu (A, \theta_i )+
\sum_{i=1}^k R_i (\mu, A, \teta_{i-1}, \theta_i, \theta_{i+1}),
\ee
with
$$ |R_i (\mu, A, \teta_{i-1}, \theta_i, \theta_{i+1}) |
\leq C_2 \exp (- C_2 L ). $$
\end{lemma}

\begin{pf}
We can write
\begin{eqnarray*}
F_\mu^k (A, \teta_1, \ldots, \teta_k )
& = & \Big( \int_{- \infty}^{\teta_1} {\cal L}_\mu (q_{A,\teta}^L (t),
{\dot{  q}}_{A,\teta}^L (t),t) +  \int_{\teta_1}^{(\teta_1 + \teta_2)/2}
{\cal L}_\mu (q_{A,\teta}^L (t), {\dot{  q}}_{A,\teta}^L (t),t) \Big) \\
& + & \sum_{i=2}^{k-1} \Big( \int_{(\teta_{i-1} + \teta_i)/2}^{\teta_i}
{\cal L}_\mu (q_{A,\teta}^L (t), {\dot q}_{A,\teta}^L (t),t) +
\int_{\teta_i}^{ ( \teta_i + \teta_{i+1} ) / 2}
{\cal L}_\mu (q_{A,\teta}^L (t), {\dot q}_{A,\teta}^L (t),t) \Big) \\
& + &
\Big( \int_{(\teta_{k-1} + \teta_k)/ 2}^{\teta_k}
{\cal L}_\mu (q_{A,\teta}^L (t), {\dot q}_{A, \teta}^L (t),t) +
\int_{\teta_k}^{+ \infty}
{\cal L}_\mu (q_{A,\teta}^L (t), {\dot q}_{A, \teta}^L (t),t) \Big).
\end{eqnarray*}
We define
$$
R_i^- (\mu, A, \teta_{i-1}, \teta_i ) =
\int_{(\teta_{i-1} + \teta_i)/2}^{\teta_i}
{\cal L}_\mu (q_{A,\teta}^L (t), {\dot q}_{A,\teta}^L (t),t) \ dt -
\int_{-\infty}^{\teta_i}
{\cal L}_\mu (q_{A,\teta_i}^\mu (t), {\dot q}_{A,\teta_i}^\mu (t),t)
\ dt,
$$
$$
R_i^+ (\mu, A, \teta_i, \teta_{i+1}) =
\int_{\teta_i}^{(\teta_i + \teta_{i+1})/2}
{\cal L}_\mu (q_{A,\teta}^L (t), {\dot q}_{A,\teta}^L (t),t) \ dt -
\int_{\teta_i}^{+ \infty}
{\cal L}_\mu (q_{A,\teta_i}^\mu (t), {\dot q}_{A,\teta_i}^\mu (t),t) \ dt
$$
where $ q_{A,\teta_i}^\mu $ is the
$ 1 $-bump pseudo-homoclinic solution obtained in lemma \ref{lem:1bump}.
Recalling the definition \ref{def:homo} of $ F_\mu (A, \teta) $ we have
\begin{eqnarray*}
F_\mu^k (A, \teta_1, \ldots, \teta_k )
& = & F_\mu (A, \teta_1) + R_1^+ ( \mu, A, \teta_1, \teta_2 ) \\
& + & \sum_{i=2}^{k-1}  F_\mu (A, \teta_i) +
\Big( R_i^- ( \mu, A, \teta_{i-1}, \teta_i)+
R_i^+ ( \mu, A, \teta_i, \teta_{i+1} \Big) \\
& + & F_\mu (A, \teta_k) + R_k^- ( \mu, A, \teta_{k-1}, \teta_k ).
\end{eqnarray*}
Setting $ R_i =   R_i^- + R_i^+ $ we derive the expression (\ref{eq:somma}).
In order to   complete the proof, it is enough  to show
the existence of $ C_2, L_2 > 0 $ such that
$ \forall \om \in {\bf R}^n $, for all $0<  \mu  \leq \mu_0 $,
$ \forall L > L_2 $, $ \forall \teta_1 < \ldots < \teta_k $ with
$ \min_i ( \teta_{i+1} - \teta_i ) > L $,
for all $ i = 1, \ldots, k $
\be\label{eq:resti}
| R_i^\pm (\mu, A, \theta_i, \theta_{i+1}) | \leq   
C_2 \exp(- C_2 L ).
\ee
We   write the proof for $ R_i^+ $.
We have
\be\label{eq:som}
\begin{array}{rcl}
R_i^+ (\mu, A, \theta_i, \theta_{i+1})& =& \dps
 \int_{\teta_i}^{(\teta_i + \teta_{i+1})/2 }
({\cal L}_\mu (q_{A,\teta}^L (t), {\dot q}_{A,\teta}^L (t),t) -
{\cal L}_\mu (q_{A,\teta_i}^{\mu} (t), {\dot q}_{A,\teta_i}^{\mu} (t),t)
 \ dt \\ &  &  \\
&-& \dps \int_{(\teta_i + \teta_{i+1})/2}^{+\infty}
{\cal L}_\mu ( q_{A,\teta_i}^\mu (t),{\dot q}_{A,\teta_i}^\mu
(t),t) \  dt .
\end{array}
\ee
  By lemma \ref{lem:1bump}-$(ii)$ the homoclinic orbit satisfies
$ \max (|q_{A,\teta_i}^\mu (t)|, |{\dot q}_{A,\teta_i}^\mu (t)|)
\leq C \exp ({- | t - \teta_i | / 2 })$.
Hence, for all $ \teta_1 < \ldots < \teta_k $ with
$ \min_i ( \teta_{i+1} - \teta_i ) > L $,
\be\label{eq:add1}
\Big| \int_{(\teta_i + \teta_{i+1})/2}^{+\infty}
{\cal L}_\mu (q_{A,\teta_i}^\mu (t), {\dot q}_{A,\teta_i}^\mu (t),t) \
dt \Big|
= O( e^{-  L/2}).
\ee
From lemma \ref{lem:kheter}-$(ii)$ we also deduce that
\be\label{eq:add2}
\Big( \int_{\teta_i}^{(\teta_i + \teta_{i+1})/2}
{\cal L}_\mu (q_{A,\teta}^L (t), {\dot q}_{A,\teta}^L (t),t)  -
{\cal L}_\mu (q_{A,\teta_i}^\mu (t), {\dot q}_{A,\teta_i}^\mu (t),t) \ dt
\Big)
= O( e^{- C L}).
\ee
From (\ref{eq:som}), (\ref{eq:add1}) and (\ref{eq:add2}) we deduce
(\ref{eq:resti}) and hence the lemma.
\end{pf}

\subsection{The diffusion orbit}

We are now able to  consider the existence of the shadowing orbit.
We give an example of  condition on $G_{\mu}$
which implies the existence of diffusion orbits.

\begin{condition}\label{spli}
({\it ``Splitting condition''})
There exist $A_0 \in {\bf T}^n $, $ \delta > 0$, $ 0 < \a < \rho $ such that
\begin{itemize}
\item $(i)$
$\inf_{\partial B_\rho (A_0)} G_\mu \geq \inf_{B_\rho ( A_0 )} G_\mu + \delta$;
\item $(ii)$
$ \sup_{ B_\alpha (A_0) } G_\mu \leq \frac{\delta}{4}
+ \inf_{B_\rho (A_0)} G_\mu $;
\item $(iii)$
$d(\{ A \in B_\rho ( A_0 ) \ | \
G_\mu (A)  \leq \delta / 2  + \inf_{B_\rho ( A_0 )} G_\mu  \},
\{ A \in B_\rho ( A_0 ) \ | \ G_\mu (A) \geq 3 \delta / 4
+  \inf_{B_\rho ( A_0 )} G_\mu \})
\geq 2 \alpha $.
\end{itemize}
\end{condition}

\begin{remark}
If $ G_\mu $ possesses a non-degenerate minimum in $ A_0 $ the ``splitting
condition'' above is satisfied, for $ \rho $ sufficiently small, 
choosing $ \d = (\min \lambda_i) \rho^2/ 4 $ and 
$ \a = ( \rho /8 ) \sqrt{(\min_i \lambda_i)/(\max_i \lambda_i) } $
where $ \lambda_i $ are the positive eigenvalues of $ D^2 G^\mu (A_0)$.
\end{remark}

\begin{remark} \label{respli}
$B_\rho (A_0)$, open ball of radius $\rho$ in   ${\bf R}^n$
(the covering space of ${\bf T}^n$),
could be replaced by a bounded open subset $U$ of ${\bf R}^n$.
\end{remark}

The following shadowing type theorem holds

\begin{theorem}\label{thm:main}
Assume $ ( H1 ) $ and the ``splitting condition'' \ref{spli}.
Then $ \forall I_0, I_0' $ with $ \om \cdot I_0= \om \cdot I_0' $,
there is a heteroclinic orbit connecting the invariant tori
${\cal T}_{I_0}$ and ${\cal T}_{I_0'}$.
Moreover there exists $ C_3 >0 $ such that
$ \forall \eta > 0 $ small enough the ``diffusion time'' $T_d$ needed
to go from a $\eta$-neighbourhood of
${\cal T}_{I_0}$ to a $\eta$-neighbourhood of
${\cal T}_{I_0'}$ is bounded by
\be\label{timediff}
T_d \leq C_3 \frac{| I_0 - I_0' |}{\delta} \rho \max \Big( | \ln \delta |,
\frac{1}{  \gamma \alpha^\tau} \Big) + C_3 | \ln ( \eta )|.
\ee
\end{theorem}

\begin{remark}
The meaning of (\ref{timediff}) is the following:
the diffusion time $ T_d $ is estimated by the product
of the number of heteroclinic transitions
$ k = $ ( heteroclinic jump / splitting ) $ = | I_0' - I_0 | / \delta $,
and of  the time $ T_s $ required for a single transition,
that is $ T_d = k \cdot T_s  $.
The time for a single transition $ T_s $ is bounded by the maximum
time between the
``ergodization time'' $ ( 1 / \gamma \alpha^\tau ) $, i.e.
the time needed for the flow $ \om t $ to make an
$ \alpha $-net of the torus,
and the time $ | \ln \delta | $ needed  to ``shadow''
homoclinic orbits for the forced pendulum equation. We use
here that these homoclinic orbits are exponentially
asymptotic to the equilibrium.
\end{remark}

\begin{remark}
The following proof works if $ G_\mu $ possesses a local 
maximum which satisfies a non-degeneracy type condition like
the ``splitting condition'' \ref{spli}, 
while in the approaches developed in \cite{BCV} and \cite{Xi},
based on Mather's theory, diffusion orbits are always built from
local minima of $ G_{\mu} $. The proof of the shadowing theorem
when the homoclimic point $ A_0 $ is a saddle point requires slightly
different arguments. For example it holds 
assuming as in \cite{G1} the condition
$ D^2 G_\mu (A_0) \om \cdot \om \neq 0 $.
\end{remark}

\begin{pf}
Assume  with no loss of generality that
$ A_0 = 0 $ and $ \inf_{ B_\rho (0) } G_\mu (A) = 0 $.
Let us choose the number of bumps $ k $ as
\be\label{eq:nbum}
k = \Big[ \frac{ 24 \cdot \rho \cdot |I_0'-I_0|}{\delta} \Big] + 1.
\ee
By lemma \ref{lem:kheter}-($i$) and  lemma \ref{lem:1bump}-$(ii)$,
the trajectory converges exponentially
fast to ${\cal T}_{I_0}$ (resp. ${\cal T}_{I_0'}$) as
$ t \to  - \infty $ (resp. $ + \infty $)
from $ \theta_1 $
(resp. $ \theta_k $). Therefore
it is enough to prove the existence of a critical point
$ ( \ov{A}, \ov{\theta}) \in {\bf T}^n \times {\bf R}^k $
of the $ k $-bump heteroclinic function $ { \cal F }_\mu^k $,
defined in (\ref{eq:kbhf}), such that for some positive constant $
  K_1 $
\be\label{totaltime}
| \ov{\teta}_k - \ov{\teta}_1 |
  \leq   K_1 \frac{| I_0 - I_0' |}{\delta} \ \rho
  \max \Big( | \ln \delta|, \frac{1}{  \gamma \alpha^\tau} \Big).
\ee
More precisely we shall   enforce
\be \label{esti}
  K_2 | \ln \delta | < | \ov{\theta}_{i+1} - \ov{\theta}_i | <
K_3
\max \Big( | \ln \delta |, \frac{1}{  \gamma \alpha^\tau} \Big) \qquad
\forall
i = 1, \ldots, k,
\ee
  for some positive constants $K_2$,$K_3$.
Let
$( \Omega_1, \cdots, \Omega_n ) $ be an orthonormal basis
of $ { \bf R }^n $ where
$$
\Omega_1 = \frac{\omega}{|\omega|}  \quad  {\rm and}  \quad
\Omega_2 = \frac{I_0'-I_0}{|I_0'-I_0|};
$$
  We recall that  $ \om \cdot (I_0'-I_0) = 0 $.
In order to find a critical point of ${\cal F}_\mu^k$
we introduce suitable coordinates
$ ( a_1, \ldots, a_n, s_1, \ldots, s_k) \in {\bf R}^n \times
(- \rho, \rho)^k $ defined by
$$
A = \sum_{j = 1}^n a_j \Omega_j,  \qquad
\teta_i = \frac{ \eta_i + s_i - a_1 }{ | \om |} \quad \forall
i = 1, \ldots, k
$$
where $ \eta_i $ are constants to be chosen later.
In these new coordinates the heteroclinic function defined in (\ref{eq:kbhf})
is given by
\be
{\wtilde {\cal F}}_\mu^k ( a_1, a_2, \ldots, a_n, s_1, \ldots, s_k) =
F_\mu \Big( \sum_{j = 1}^n a_j \Omega_j ,
\frac{ \eta_1 + s_1 - a_1}{| \omega |}, \ldots,
\frac{ \eta_k + s_k - a_1}{| \omega |} \Big) - |I_0'-I_0| a_2.
\ee
Using the invariance property (\ref{invariance}) we see that
${\wtilde {\cal F}}_\mu^k $ does not depend on the new
variable $ a_1 $   :
\begin{eqnarray*}
{\wtilde {\cal F}}_\mu^k (a_1, a_2, \ldots, a_n, s_1, \ldots, s_k) & = &
F_\mu \Big( \sum_{j = 2}^n a_j \Omega_j ,
\frac{ \eta_1 + s_1}{| \omega |}, \ldots,
\frac{ \eta_k + s_k}{| \omega |} \Big) - |I_0'-I_0| a_2 \\
& = &
{\wtilde {\cal F}}_\mu^k (0, a_2, \ldots, a_n, s_1, \ldots, s_k).
\end{eqnarray*}
For simplicity of notation we will still denote
$
{\wtilde {\cal F}}_\mu^k (a_2, \ldots, a_n, s_1, \ldots, s_k) :=
{\wtilde {\cal F}}_\mu^k (0, a_2, \ldots, a_n, s_1, \ldots, s_k).
$
We now choose the contants $ \eta_i $. Let
\be\label{eq:sep}
D =  \frac{ | \om | }{C_2} \Big| \ln \Big(\frac{24 C_2}{\delta} \Big)
\Big| + 2 \rho,
\ee
where $C_2$ is the constant appearing in lemma \ref{approxsum}.
We shall use the following fact (see \cite{BGW}):  there is
$ \ov{C} >0 $ such that, for all intervals $ J \subset {\bf R} $ of
length greater or equal to $  \ov{C} / (\gamma \alpha^\tau) $, there is
$ \theta \in J $ such that
\be \label{number}
d( \theta \omega, 2 \pi {\bf Z}^n) < \alpha.
\ee
By (\ref{number}) there   is  $ (\eta_1, \ldots, \eta_k ) \in {\bf R}^k $
such that
\be\label{theta2}
\eta_i \Omega_1 \equiv \chi_i, \ {\rm mod}  2 \pi {\bf Z}^n,
\quad | \chi_i | < \alpha \quad {\rm and } \quad \chi_i \cdot \Omega_1 = 0,
\ i.e. \ \chi_i = \sum_{j = 2}^n \chi_{i,j} \Omega_j.
\ee
\be \label{theta1}
\eta_1 = 0,  \  \  \  D \leq \eta_{i+1}- \eta_i
\leq \Big(D  +   \frac{ \ov{C} |\om |}{\gamma \alpha^{\tau} } \Big).
\ee
By (\ref{eq:sep}), (\ref{theta1}), since $s_i \in (- \rho, \rho)$
we have that $ \teta_{i+1} - \teta_i \geq
\frac{ 1 }{C_2} | \ln (\frac{24 C_2}{\delta})| $;
hence, by lemma \ref{approxsum}, setting
\begin{eqnarray*}
{\wtilde R}_i & = & {\wtilde R}_i (a_2, \ldots, a_n, s_{i-1}, s_i, s_{i+1} )
\\
& = & R_i \Big( \sum_{j = 2}^n a_j \Omega_j,
\frac{s_{i-1} + \eta_{i-1} - a_1}{ |\om |},
\frac{s_i + \eta_i - a_1}{|\om |},
\frac{s_{i+1} + \eta_{i+1} - a_1}{ |\om |} \Big)
\end{eqnarray*}
we   get
\be\label{eq:smr}
|{\wtilde R}_i
(a_2, \ldots, a_n, s_{i - 1}, s_i, s_{i+1}) | \leq \frac{ \delta}{24}.
\ee
By lemma \ref{approxsum}, the invariance property (\ref{invariance}),
(\ref{theta2}) and since
$ G_\mu $   is $2\pi {\bf Z}^n$-periodic, we have
\begin{eqnarray*}
{\wtilde {\cal F}}_\mu^k (a_2, \ldots ,a_n, s_1, \ldots, s_k ) & = &
\sum_{i = 1}^k F_\mu \Big( \sum_{j = 2}^n a_j \Omega_j,
\frac{ \eta_i + s_i }{ | \omega |} \Big) +
{\wtilde R}_i - |I_0'-I_0| a_2 \\
& = & \sum_{i = 1}^k F_\mu \Big( \sum_{j = 2}^n a_j \Omega_j +
\chi_i + s_i \Omega_1, 0   \Big) +
{\wtilde R}_i - |I_0'-I_0| a_2 \\
& = & \sum_{i = 1}^k G_\mu \Big( \sum_{j = 2}^n
( a_j + \chi_{i,j} ) \Omega_j + s_i \Omega_1 \Big) +
{ \wtilde R }_i - |I_0'-I_0| a_2 \\
& = &
\sum_{i = 1}^k
{\wtilde G}_\mu ( a_2 + \chi_{i,2}, \ldots, a_n + \chi_{i,n}, s_i ) +
{\wtilde R}_i - |I_0'-I_0| a_2
\end{eqnarray*}
where
$ {\wtilde G}_\mu ( a_2, \ldots, a_n, s) =
G_\mu ( \sum_{j = 2}^n a_j \Omega_j + s \Omega_1 ) $.
Since the basis $ ( \Omega_1, \ldots, \Omega_n) $ is orthonormal
the function $ { \wtilde G}_{\mu} $
satisfies the same properties   as  $ G_\mu $, {\it i.e.}
$$
\sup_{ B_\alpha (0)}{\wtilde G }_\mu \leq \delta / 4, \
\inf_{ \partial B_\rho (0)} \wG_{\mu} \geq \delta \ {\rm and} \
d(\{ x \in  B_\rho (0) \ | \
{\wtilde G}_\mu (x) \leq \delta / 2 \}, \{ x \in B_\rho (0) \ | \
{\wtilde G}_\mu (x) \geq 3 \delta / 4 \}) \geq 2 \alpha.
$$
We shall find a critical point of $ {\wtilde {\cal F}}_\mu^k $ in
$$
W = \Big\{ ( a_2, \ldots, a_n, s)
\in {\bf R}^{ n - 1 } \times {\bf R}^k \ \Big|
\  ( a_2, \ldots, a_n, s_i) \in B_\rho ( 0 ),
\quad \forall i =1, \ldots, k \Big\}.
$$
$ {\wtilde {\cal F}}_\mu^k $ attains its minimum over $ \ov{W} $ at some
point
$ ( \ov{a}, \ov{s})$. Notice that by (\ref{eq:smr})
$$
\inf_{\ov{W}} {\wtilde {\cal F}}_\mu^k \leq {\wtilde {\cal F}}_\mu^k ( 0 ,
0 ) =
\sum_{i=1}^k \wG_{\mu} (0,\chi_i) + k \frac{\delta}{24}.
$$
  Since
$ | \chi_i | < \alpha $ for all $ i = 1, \ldots, k $
and $\sup_{ B_\alpha (0)}{\wtilde G }_\mu \leq \delta / 4$, we have
\be \label{inf}
\inf_{ \ov { W }} {\wtilde {\cal F}}_\mu^k \leq k\frac{\delta}{4} +
k \frac{\delta}{24} = k \frac{7 \delta}{24}.
\ee
The theorem is proved if we show that
$(\ov{a}, \ov{s}) \in W $.
Arguing by contradiction assume that $(\ov{a}, \ov{s}) \in
\partial W$.   Then there is some $ l \in \{1, \cdots, k  \}$ such that
$ (\ov{a} + \chi_l, \ov{s}_l) \in \partial B_\rho (0)$,
  so that
$\wG_{\mu}(\ov{a} + \chi_l, \ov{s}_l) \geq \delta $.
We now prove that
$( \ov{a} + \chi_l, t ) ; t \in (- \rho, \rho ) \}
\cap B_{\rho} (0)  \subset  Z := \{ x \in B_\rho (0) \ | \
\wG_{\mu} ( x )   \geq 3 \delta/4  \}$. Indeed, if not, by
  (\ref{eq:smr}), for some  $t\in (-\rho,\rho)$ such that
$( \ov{a} + \chi_l, t ) \in B_{\rho} (0)$,
\begin{eqnarray*}
{\wtilde {\cal F}}_\mu^k( \ov{a},
\ov{s}_1, \cdots , \ov{s}_{l-1},
t, \ov{s}_{l+1}, \cdots , \ov{s}_k) & \leq &
{\wtilde {\cal F}}_\mu^k(\ov{a},\ov{s})+
(\wG_{\mu}(\ov{a}+ \chi_l,t) -
\wG_{\mu}(\ov{a}+ \chi_l, \ov{s}_l)))\\
& +& |{\wtilde R}_{l-1} (\ov{a}, \ov{s}_{l-2},
\ov{s}_{l-1}, s_l) -
{\wtilde R}_{l-1}
(\ov{a},\ov{s}_{l-2},\ov{s}_{l-1},t)|\\
& + & | {\wtilde R}_l (\ov{a}, s_{l-1}, t,\ov{s}_l ) -
{\wtilde R}_l ( \ov{a}, \ov{s}_{l-1},
\ov{s}_l,\ov{s}_{l+1})| \\
& + & | {\wtilde R}_{l+1} (\ov{a}, \ov{s}_l,
\ov{s}_{l+1}, t ) -
{\wtilde R}_{l+1} ( \ov{a}, \ov{s}_l,
\ov{s}_{l+1},\ov{s}_{l+2})| \\
& < &
{\wtilde {\cal F}}_\mu^k( \ov{a}, \ov{s}) - \frac{\delta}{4} +
\frac{6 \delta }{24} = {\wtilde {\cal
F}}_\mu^k(\ov{a},\ov{s}),
\end{eqnarray*}
which is   wrong since $(\ov{a}, \ov{s})$ is the minimum
of ${\wtilde {\cal F}}_\mu^k$ over $ \ov{W} $. We deduce in   particular
that, for all $i$,
$(\ov{a} + \chi_l, \ov{s_i}) \in Z \cup B_{\rho}(0)^c$.
  Now, as $\wtilde{G}_{\mu} \geq 3\delta/4$ in a neighbourhood
of $\partial B_{\rho}(0)$, our splitting condition implies that
$$
  d ( \{ x \in B_\rho (0) \ | \
\wG_{\mu} ( x )  \geq  3 \delta/4 \} \cup B_{\rho}(0)^c , \{ x \in B_\rho(0)
\ | \
\wG_{\mu} (x) \leq \delta/2 \}) \geq  2\alpha.
$$
We derive   by (\ref{theta2}) that for all $i$,
\be\label{eq:supper}
\wG_{\mu}(\ov{a} + \chi_i , \ov{s_i}) \geq
\delta/2.
\ee
As a consequence, noting that, from (\ref{eq:nbum}),
$ | I_0' - I_0 | \ \rho \leq (k \delta )/ 24 $, we deduce that
$$
{\wtilde {\cal F}}_\mu^k(\ov{a}, \ov{s}) \geq k
\frac{\delta}{2} -
k\frac{\delta}{24} - k \frac{\delta}{24} = k \frac{10 \delta}{24} >
k \frac{7 \delta}{24},
$$
contradicting (\ref{inf}).
The proof of the theorem is complete.
\end{pf}

When the frequency vector  $\om $ is considered as a
constant, independent of any parameter (``a priori-unstable case'')
it is easy to justify the splitting condition \ref{spli} using the
first-order approximation given by the Poincar\'e-Melnikov primitive.
With a Taylor expansion in $ \mu $ we can easily prove that
for $ \mu $ small enough
$$
G_\mu (B) = Const + \mu \Gamma (B) + O(\mu^2), \ \forall B \in {\bf T}^n,
$$
where $ \Gamma: {\bf T}^n \to {\bf R} $
is nothing but the Poincar\'e-Melnikov primitive
$$
\Gamma (B) =
\int_{{\bf R}} (1- \cos q_0 (t) ) f (\om t + B) \ dt.
$$
\noindent
Hence, if $ \Gamma $ possesses a proper minimum (resp. maximum) 
in $ A_0 \in {\bf R}^n $, i.e  $ \exists  r > 0 $ such that   
$ \inf_{\partial B_r ( A_0 )} \Gamma > \Gamma (A_0 ) $
(resp. $ \sup_{\partial B_r ( A_0 )} \Gamma < \Gamma ( A_0 ) $) 
then, for $ \mu $ small enough,
the ``splitting'' condition \ref{spli} holds with
$ \delta = O(\mu) $, $ \rho = O(1) $ and
$ \alpha = O( 1) $. We remark that the previous 
$ B_r (A_0 ) $ could be replaced 
by a bounded open subset $ U $ of $ { \bf R }^n $. 
Applying theorem \ref{thm:main} we deduce

\begin{theorem}\label{thm:apu}
Assume $ (H1) $ and let $ \Gamma $ possess a proper minimum 
(or maximum) $ A_0 $, i.e. $ \exists  r > 0 $ such that   
$ \inf_{\partial B_r ( A_0 )} \Gamma > \Gamma ( A_0 ) $. Then, 
for $ \mu $ small enough,
the same statement of theorem \ref{thm:main} holds with a diffusion time
$ T_d = O( ( 1 / \mu ) \log (1 / \mu )) $.
\end{theorem}

\section{More general perturbation terms}

In this section we show how to adapt the arguments of the
previous section when dealing with a more general perturbation
term  $f( \vphi, q)$. Regarding regularity it is sufficient
to have finite large enough smoothness for $ f $.
The equation of motion derived by Hamiltonian ${\cal H}_\mu$ are
\be\label{eqmotion1}
\dot{\vphi} = \om, \qquad
\dot{I}  = - \mu \partial_{\vphi} f ( \vphi, q), \qquad
\dot{q} = p, \qquad
\dot{p} = \sin{q} - \mu \partial_q f ( \vphi, q),
\ee
corresponding to the quasi-periodically forced pendulum
\be\label{forcedpend}
- \ddot{q} + \sin{q} = \mu \ \partial_q f ( \om t + A, q).
\ee

\subsection{Invariant tori in the perturbed system}

The first step is to prove the persistence of invariant tori
for $ \mu \neq 0 $ small enough.
It appears that no more than the standard Implicit Function Theorem
is required to prove the following well known result
(see for example \cite{Gen} for a proof)

\begin{theorem}\label{thm:invto}
Let $ \om $ satisfy $(H1)$. For $\mu$ small enough
and $ \forall I_0 \in {\bf R}^n $
system ${\cal H}_\mu $ possesses $ n $-dimensional invariant tori
$ {\cal T}_{I_0}^\mu \approx {\cal T}_{I_0} $ of the form
\be\label{toripert}
{\cal T}_{I_0}^\mu = \Big\{
I = I_0 + a^\mu ( \psi ), \ \vphi = \psi, \
q =  Q^\mu ( \psi ), \
p =  P^\mu ( \psi ), \quad \psi \in {\bf T}^n \Big\},
\ee
with $ Q^\mu (\cdot ), P^\mu (\cdot)  = O( \mu ) $,
$ a^\mu (\psi ) = O( \mu )$.
Moreover the dynamics on $ {\cal T}_{I_0}^\mu $ is
conjugated to the rotation of speed $ \om $ for $\psi$.
\end{theorem}

We first determine the functions
$ Q^\mu ( \cdot ), P^\mu ( \cdot )$ in (\ref{toripert}).
Using the standard Implicit Function Theorem
we prove that there exists a unique quasi-periodic solution $q^\mu_{A} (t)$
for the quasi-periodically forced pendulum (\ref{forcedpend})
which bifurcates from the hyperbolic equilibrium $0$.

\begin{lemma}
Let $f \in C^l ({\bf T}^n \times {\bf T}) $.
For $ \mu $ small enough there exists a unique quasi-periodic solution
$(q_A^\mu (t), p_A^\mu (t))$ of (\ref{forcedpend})
with $(q_A^\mu (t), p_A^\mu (t)) = O( \mu ) $, $C^{l-1}$-smooth
in $ A $. More precisely there exist
functions $ Q^\mu, P^\mu: {\bf T}^n \to {\bf R}$ of class $C^{l-1}$,
such that
$(q_A^\mu (t), p_A^\mu (t))= (Q^\mu (\omega t + A) , P^\mu (\omega t + A)).$
\end{lemma}

\begin{pf}
Let $L$ be the Green operator of the
differential operator $ h \to - D^2 h  + h $
with Dirichlet boundary conditions at $\pm \infty$.
$L$ is explicitely given by
$ L (f) = \int_{{\bf R}} e^{-|t-s|} f(s) \ ds / 2$.
It results that $L$ is a continuous linear operator in
the Banach space of the continuous bounded functions
from $\bf{R}$ to $\bf{R}$, which we shall denote
by $E$. We consider the non-linear operator
$ S: {\bf R} \times {\bf T}^n \times
E \to  E $
$$
S (\mu, A, q) := q - L( q- \sin q) - \mu L( \partial_q f (\om t + A, q)).
$$
$S$ is of class $C^{l-1}$. We are looking  for a solution
$ q^\mu_A $ of $ S (\mu,A, q)=0 $. Since $ S (0,A, 0)=0 $ and
$ \partial_q S (0,A, 0) = Id $,
by the Implicit Function Theorem there exists, for $ \mu $ small enough,
a unique solution  $ q^{\mu}_A = O(\mu) $. By (\ref{forcedpend})
$q_A^{\mu} \in C^{l+1}({\bf R})$; moreover it is $C^{l-1}$-smooth
in $ A $. We define the $C^{l-1}$-maps
$ Q^\mu ( \cdot ), P^\mu ( \cdot ) : {\bf T}^n \to {\bf R}$ by
$$
Q^\mu ( A  ) :=  q^\mu_A (0),
\qquad P^\mu ( A ) :=  \dot{q}^\mu_A (0).
$$
By uniqueness we deduce that
$ q^\mu_{ A} ( s + t ) = q^\mu_{A + \omega s} (t), \ \forall s, t \in {\bf
R}$.
For $ t =  0 $ this yields
$$
q^\mu_{A} ( s ) =
q^\mu_{A + \omega s } (0) := Q^\mu ( A + \omega s) \quad {\rm and} \quad
p^\mu_{A} ( s ) =
p^\mu_{A + \om s} (0) := P^\mu ( A + \om s), \ \forall  s \in {\bf R}
$$
proving the lemma.\end{pf}

We now define the functions $a^\mu (\psi )$ of (\ref{toripert}).
We impose that $(\om t+A, I_0 + a^\mu ( \om t + A ), Q^{\mu} (\om
t+A), P^{\mu} (\om t+A))$
satisfy the equations of motions (\ref{eqmotion1}); hence the functions
$ a^\mu ( \psi ) $ must satisfy the following system of equations
\be\label{eq:diof}
(\om \cdot \nabla) a^\mu ( \psi ) = \mu g^\mu ( \psi ), \quad {\rm where}
\quad g^\mu ( \psi ) := - (\nabla_{\psi} f) ( \psi, Q^\mu ( \psi)).
\ee
In order to solve (\ref{eq:diof}) we expand in Fourier series
the functions
$ a^\mu ( \psi ) = \sum_{k \in {\bf Z}^n} a_k e^{i k \cdot \psi}$,
$ g^\mu ( \psi ) = \sum_{k \in {\bf Z}^n} g_k e^{i k \cdot \psi}$.
Each Fourier coefficient $ a_k $ must then satisfy
\be\label{eq:fourier}
i (k \cdot \om )  a_k = \mu g_k,  \quad \quad \forall  k \in {\bf Z}^n.
\ee
It is necessary for the existence of a solution that 
$g_0 = \int_{{\bf T}^n}  g^\mu (\psi ) d \psi =0 $.
This property can be checked directly, that is

\begin{lemma}
We have
\be\label{eq:zeromean}
\int_{{\bf T}^n}  (\nabla_{\psi} f)
(  \psi, Q^\mu ( \psi )) \ d \psi = 0.
\ee
\end{lemma}

\begin{pf}
For all $i = 1, \ldots, n $
\be\label{eq:firstco}
\partial_{\psi_i}
f (\psi, Q^\mu ( \psi)) =
\frac{d}{d \psi_i }
f (\psi, Q^\mu ( \psi) ) -
\partial_q f (\psi, Q^\mu ( \psi))
\partial_{\psi_i} Q^\mu ( \psi ).
\ee
Since $( q^\mu_A (t), p^\mu_A (t) )$ satisfies the pendulum equation
$ \sum_{j=1}^n \omega_j \partial_{\psi_j} Q^\mu (\psi)  = P^\mu (\psi)$,
$ \sum_{j=1}^n \omega_j \partial_{\psi_j} P^\mu (\psi)  = \sin Q^\mu
\psi)  -
\mu \partial_q f ( \psi,  Q^\mu (\psi))$ and we deduce that
\be\label{eq:part1}
- \partial_q f (\psi, Q^\mu ( \psi ))
\partial_{\psi_i} Q^\mu ( \psi ) =
\frac{1}{\mu} \Big( \frac{d}{d \psi_i} \cos Q^\mu (\psi ) +
\sum_{j=1}^n \omega_j \partial_{\psi_j} P^\mu \partial_{\psi_i} Q^\mu
(\psi )
\Big).
\ee
We now prove that
\be\label{eq:pez1}
\sum_{j = 1}^n \omega_j \partial_{\psi_j} P^\mu \partial_{\psi_i} Q^\mu =
\frac{d}{d \psi_i} \frac{( P^\mu (\psi ))^2}{2} +
\sum_{j \neq i} \omega_j \Big( \partial_{\psi_i}
(Q^\mu \partial_{\psi_j} P^\mu ) -  \partial_{\psi_j}
(Q^\mu \partial_{\psi_i} P^\mu ) \Big).
\ee
Indeed
\begin{eqnarray*}
\sum_{j=1}^n \omega_j \partial_{\psi_j} P^\mu \partial_{\psi_i} Q^\mu & = &
\omega_i \partial_{\psi_i} P^\mu \partial_{\psi_i} Q^\mu +
\sum_{j \neq i}  \omega_j \partial_{\psi_j} P^\mu \partial_{\psi_i} Q^\mu \\
& = &
\Big( \sum_{j=1}^n \omega_j \partial_{\psi_j} Q^\mu \Big)
\partial_{\psi_i} P^\mu + \sum_{j \neq i} \omega_j
\Big( \partial_{\psi_j} P^\mu \partial_{\psi_i} Q^\mu -
\partial_{\psi_j} Q^\mu \partial_{\psi_i} P^\mu \Big) \\
& = &
\frac{d}{d \psi_i} \frac{(P^\mu (\psi))^2}{2} +
\sum_{j \neq i}  \omega_j \Big(  \partial_{\psi_i}
(Q^\mu \partial_{\psi_j} P^\mu ) -  \partial_{\psi_j}
(Q^\mu \partial_{\psi_i} P^\mu ) \Big).
\end{eqnarray*}
From (\ref{eq:firstco}), (\ref{eq:part1}), (\ref{eq:pez1}) we finally
obtain that
$$
\partial_{\psi_i} f (\psi, Q^\mu ) =
\frac{d}{d \psi_i } \Big( f (\psi, Q^\mu ) +
\frac{1}{\mu} \cos Q^\mu + \frac{1}{\mu}  \frac{{P^\mu}^2}{2} \Big)
+
\frac{1}{\mu} \sum_{j \neq i}  \omega_j \Big(  \partial_{\psi_i}
(Q^\mu \partial_{\psi_j} P^\mu ) -  \partial_{\psi_j}
(Q^\mu \partial_{\psi_i} P^\mu ) \Big)
$$
from which property (\ref{eq:zeromean}) follows.
\end{pf}

Since $ \om $ satisfies $ (H1) $ and $f $ is sufficiently smooth
the function $ a^\mu $ defined by
\be\label{eq:fouriersol}
a^\mu (\psi ) = \sum_{k \in {\bf Z}^n, k \neq 0}
\frac{g_k}{i (k \cdot \om )} e^{i k \cdot \psi},
\ee
which formally solves equation (\ref{eq:fourier}),
is well defined and smooth.
Indeed since $ f \in C^l $ 
the function $ g^\mu $ defined in (\ref{eq:diof}) is
$ C^{l-1} $ and there exists $ M > 0 $ such that
$ | g_k | \leq M / |k|^{l-1}$, $ \forall  k \in {\bf Z}^n $,
$ k \neq 0 $. By $(H1)$ it follows that
$| a_k | \leq M / |k|^{l-1} | \om \cdot k |
         \leq M |k|^{\tau} / (\gamma |k|^{l-1}).$
The proof of theorem \ref{thm:invto} is complete.

\subsection{The new symplectic coordinates}

In order to reduce to the previous case we want to put the tori
${\cal T}_{I_0}^\mu $ at the origin by a symplectic
change of variables.
Recalling that the tori ${\cal T}_{I_0}^\mu$ are
{\it isotropic} submanifolds we can prove the following lemma

\begin{lemma}
The transformation of coordinates $(J, \psi, u,v) \to (I, \vphi, q, p)$
defined on the covering space $ {\bf R}^{2 ( n + 1 )} $ of
$ {\bf T}^n \times {\bf R}^n \times {\bf T} \times {\bf R} $ by
\be \label{trascano}
I =  a^\mu (\psi ) + u  \partial_\psi P^\mu (\psi)  -
v \partial_\psi Q^\mu (\psi)  + J, \quad \vphi = \psi,
\quad  q  = Q^\mu (\psi) + u,
\quad  p  = P^\mu (\psi) + v
\ee
is symplectic.
\end{lemma}

\begin{pf}
Set $ d I \wedge d \vphi = \sum_{i=1}^n d I_i \wedge d \vphi_i $
and $d J \wedge d \psi = \sum_{i=1}^n d J_i \wedge d \psi_i$. We have
\begin{eqnarray*}
d I \wedge d \vphi + dp \wedge dq & = &
\sum_{i=1}^n d a^\mu_i (\psi) \wedge d\psi_i +
d ( u \partial_{\psi_i} P^\mu (\psi)   ) \wedge d\psi_i
- d (v  \partial_{\psi_i} Q^\mu (\psi) ) \wedge d\psi_i \\
& + &
d J \wedge d \psi + dv \wedge du +
d P^\mu (\psi) \wedge d Q^\mu (\psi) +
d P^\mu (\psi) \wedge du +  dv \wedge  d Q^\mu (\psi).
\end{eqnarray*}
Using that the tori  ${\cal T}_{I_0}^\mu $ are isotropic, that is
$ \sum_{i=1}^n d a^\mu_i (\psi) \wedge d\psi_i +
d P^\mu (\psi) \wedge  d Q^\mu (\psi) = 0$,
and noticing that
$ \sum_{i,j} u\partial_{\psi_i,\psi_j}^2
P^\mu (\psi) d\psi_j \wedge d \psi_i = 0$
$=\sum_{i,j} v\partial_{\psi_i,\psi_j}^2
Q^\mu (\psi) d\psi_j \wedge d \psi_i $
we deduce
\begin{eqnarray*}
d I \wedge d \vphi + dp \wedge dq & = &
\sum_{i=1}^n  d (u \partial_{\psi_i} P^\mu ( \psi )   ) \wedge d\psi_i
- d (v \partial_{\psi_i} Q^\mu (\psi) ) \wedge d\psi_i \\
& + &
d J \wedge d \psi + dv \wedge du +
d P^\mu (\psi) \wedge du +  dv \wedge  d Q^\mu (\psi) \\
& = & d J \wedge d \psi + dv \wedge du +
\sum_{i=1}^n \partial_{\psi_i} P^\mu (\psi) du \wedge d\psi_i -
\partial_{\psi_i} Q^\mu (\psi) dv  \wedge d \psi_i \\
&+&  d P^\mu (\psi) \wedge du +  dv \wedge  d Q^\mu (\psi) \\
& = & d J \wedge d \psi +  dv \wedge du,
\end{eqnarray*}
and the transformation (\ref{trascano}) is symplectic.
\end{pf}

In the new coordinates each invariant torus ${\cal T}_{I_0}^\mu $
is simply described by
$\{ J = I_0 , \ \psi \in {\bf T}^n, \  u = v = 0 \}$ and the
new Hamiltonian writes
$$
{\cal K}_\mu = E_\mu + \om \cdot J + \frac{v^2}{2} + ( \cos u - 1)  +
P_0 (\mu, u, \psi) \leqno{({\cal K}_{\mu})}
$$
where
$$
P_0 (\mu,  u, \psi ) = \Big( \cos( Q^\mu + u) - \cos Q^\mu +
(\sin Q^\mu ) u + 1 - \cos u \Big) +
\mu \Big( f( \psi, Q^\mu  + u) -  f ( \psi, Q^\mu ) -
\partial_q f ( \psi, Q^\mu   ) u \Big)
$$
and  $ E_\mu $ is the energy of the perturbed
invariant torus ${\cal T}_0^{\mu}=\{ ( a^\mu ( \psi ), \psi, Q^\mu (\psi
), P^\mu (\psi )); \psi \in {\bf T}^n \}.$
Hamiltonian $({\cal K}_\mu)$  corresponds to
the quasi-periodically forced pendulum equation
\be\label{pendtrasf}
- \ddot{u} + \sin{u} =  \partial_u P_0 (\mu, u,  \om t + A).
\ee
of Lagrangian
\be\label{eq:lagrap}
L_\mu =
\frac{\dot{u}^2}{2} + ( 1- \cos u) - P_0 (\mu, u, \om t + A).
\ee

Since the Hamiltonian $ {\cal K}_\mu $ is no more
periodic in the variable $ u $ we can not directly apply
theorem \ref{thm:main} and the
arguments of the previous sections require some modifications.
Arguing as in lemma \ref{lem:1bump} we deduce that,
there exists, for $ \mu $ small enough,
a unique 1-bump pseudo-homoclinic solution
$ u_{A, \teta}^\mu (t) $, true solution of (\ref{pendtrasf})
in $( - \infty, \teta ), (\teta, + \infty) $,
satisfying all the properties of lemma \ref{lem:1bump}.
Then we define the function
${\cal F}_\mu: {\bf T}^n \times {\bf R} \to {\bf R} $ as
\begin{eqnarray*}
{\cal F}_\mu ( A, \theta) & = & \int_{-\infty}^\theta
\frac{({\dot u}^\mu_{A, \theta})^2}{2} + ( 1 - \cos u^\mu_{A, \theta} )
- P_0 (\mu, u^\mu_{A, \theta} , \omega t + A) \ dt \\
& + & \int_\theta^{+\infty}
\frac{({\dot u}^\mu_{A, \theta})^2}{2} + ( 1 - \cos u^\mu_{A, \theta} )
- P_1 ( \mu, u^\mu_{A, \theta}, \om t + A) \ dt
+ 2 \pi {\dot q}^\mu_A ( \theta ),
\end{eqnarray*}
where, $ \forall i \in {\bf Z} $, we have set
\begin{eqnarray*}
P_i (\mu, u, \omega t + A) & = & \Big( \cos (q^\mu_A (t)+ u ) -
\cos q^\mu_A (t)  + \sin q^\mu_A (t) \ ( u - 2 \pi i ) + 1 - \cos u \Big) \\
& + &
\mu \Big( f(\omega t + A, q^\mu_A (t) + u) - f( \omega t + A, q^\mu_A
(t) ) -
(\partial_q f) (\omega t + A, q^\mu_A (t) ) \ (u - 2 \pi i ) \Big).
\end{eqnarray*}
Since $ u^\mu_{A, \theta} $ converges exponentially fast to $0$
for $ t \to - \infty $ and to $ 2 \pi $ for $ t \to + \infty $
the above integrals are convergent.
The term $ 2 \pi {\dot q}^\mu_A ( \theta ) $ takes into account
that the stable and the
unstable manifolds of the tori $ {\cal T}^\mu_{I_0} $ are not exact
Lagrangian manifolds, see \cite{LMS}.
We define the {\it ``homoclinic function''}
${\cal G}_\mu: {\bf T}^n \to {\bf R}$
as
\be\label{eq:homospligp}
{\cal G}_\mu ( A ) = {\cal F}_\mu ( A, 0 ).
\ee
It holds also
$ {\cal F}_\mu  (A , \teta ) = {\cal G}_\mu (A + \om \teta), \quad \forall
\teta \in {\bf R}.$
Arguing as in lemma \ref{lem:kheter} we can prove the existence
of $ k $-bump pseudo-homoclinic solutions $ u_{ A, \teta }^L $,
which is a true solution of (\ref{pendtrasf})
in each interval $ ( - \infty, \teta_1 ) $, $ ( \teta_i, \teta_{i+1})$
($ i = 1, \ldots, k-1 $), $ ( \teta_k, +\infty ) $,
and satisfying all the properties of lemma \ref{lem:kheter}.
Then we define the ``$ k $-bump heteroclinic function''
\begin{eqnarray*}
{\cal F}_\mu^k ( A , \teta_1, \ldots, \teta_k ) & = &
\int_{-\infty}^{\theta_1}
\frac{({\dot u}^L_{A, \theta})^2}{2} + ( 1 - \cos u^L_{A, \theta} )
- P_0 (\mu, u^L_{A, \theta} , \omega t + A) \ dt  +
2 \pi {\dot q}^\mu_A ( \teta_1 ) \\
& + & \sum_{i = 1}^{k-1} \int_{\teta_i}^{\teta_{i+1}}
\frac{({\dot u}^L_{A, \theta})^2}{2} + ( 1 - \cos u^L_{A, \theta} )
- P_i (\mu, u^L_{A, \theta} , \omega t + A) \ dt +
2 \pi {\dot q}^\mu_A ( \teta_{i+1} ) \\
& + & \int_{\theta_k}^{+\infty}
\frac{({\dot u}^L_{A, \theta})^2}{2} + ( 1 - \cos u^L_{A, \theta} )
- P_k ( \mu, u^\mu_{A, \theta}, \om t + A) \ dt - (I_0'- I_0) \cdot A
\end{eqnarray*}

If $ \partial_{\teta_i}{\cal F}_\mu^k (A, \teta_1, \ldots ,\teta_k ) =
({\dot u}^L_{A, \teta})^2 ( \teta_i^- ) / 2 -
({\dot u}^L_{A, \teta})^2 ( \teta_i^+ ) / 2 = 0 $
then $ u^L_{A, \teta} $ is a true solution
of the quasi-periodically forced pendulum (\ref{pendtrasf}).
As in the previous section
the variation in the action variables is given by the partial
derivative with respect to $ A $, that is
\be\label{jumpgt}
\partial_A {\cal F}_\mu^k ( A, \teta )= \int_{-\infty}^{+\infty}
- \mu \Big( \partial_\varphi f (\omega t + A, q^\mu_A (t) +
u^L_{A, \theta} (t) ) -
\partial_\varphi f (\omega t + A, q^\mu_A (t)) \Big) \ dt -
(I_0'- I_0).
\ee

\begin{lemma}\label{lem:vnk}
Let $ ( A , \teta ) $ be a critical point of
$ {\cal F}_\mu^k $. Then there exists
a heteroclinic orbit connecting the tori
${\cal T}^\mu_{I_0}$ and ${\cal T}^\mu_{I_0'}$.
\end{lemma}

\begin{pf}
By (\ref{jumpgt}) it is easy to verify
that the solutions of (\ref{eqmotion1})
$ ( I_\mu (t), \om t + A, q^\mu_A +  u^L_{A, \theta},
\dot{q}^\mu_A +  {\dot u}^L_{A, \theta} )$,
with
$ I_\mu (t) = C - \mu \int_0^t \partial_\vphi f (\om s + A,
   q^\mu_A (s) +  u^L_{A, \theta} (s)) \ ds$
and
$ C = I_0' + a_\mu (A) + \mu \int_0^{+ \infty}
\partial_\vphi f (\om t + A, q^\mu_A (t) + u^L_{A, \theta}(t)) -
\partial_\vphi f (\om t + A, q^\mu_A (t) ) \ dt, $
is a heteroclinic solution connecting
${\cal T}^\mu_{I_0}$ and ${\cal T}^\mu_{I_0'}$.
\end{pf}

Finally, arguing as in the proof of theorem \ref{thm:main}, we obtain

\begin{theorem}\label{thm:main1g}
Assume $ ( H1 ) $ and let ${\cal G}_\mu $ satisfy the ``splitting
condition'' \ref{spli}.
Then $\forall I_0, I_0'$ with $\om \cdot I_0= \om \cdot I_0'$,
there is a heteroclinic orbit connecting the invariant tori
${\cal T}_{I_0}^\mu $ and ${\cal T}_{I_0'}^\mu $.
The same estimate on the diffusion time
given in theorem \ref{thm:main} holds.
\end{theorem}

A Taylor expansion in $\mu $ gives

\begin{lemma}\label{lem:pmg}
For $ \mu $ small enough
\be\label{eq:eqsvi}
{\cal G}_\mu (A ) = const + \mu M (A) + O(\mu^2), \qquad \forall A \in
{\bf T}^n
\ee
where
$ M ( A ) $ is the Poincar\'e-Melnikov primitive
$
M ( A ) = \int_{-\infty}^{+\infty}
\Big[ f ( \omega t + A, q_0 (t)) - f ( \omega t + A, 0 ) \Big] \ dt.
$
\end{lemma}

\begin{pf}
We develop with a Taylor expansion in $ \mu $ the
Lagrangian $ L _\mu $ defined in (\ref{eq:lagrap})
\be\label{svilagra}
L_\mu = \frac{\dot{u}^2}{2} + ( 1- \cos u) +
\mu \Big( (u - \sin u ) \g +
f ( \omega t + A, u ) - f ( \omega t + A, 0 ) -
\partial_q f (\omega t + A , 0 )u
\Big) + {\cal R}(\mu, u, t)
\ee
where $ \g (t) := {\partial_\mu}_{| \mu = 0} q^\mu_A (t)$,
$ |{\cal R}(\mu, u, t)| = (\mu^2 )$,
$ {\cal R}(\mu, 0, t) =0 $ and $\partial_u {\cal R}(\mu, 0, t) =0 $.
The Melnikov function corresponding to Lagrangian
(\ref{svilagra}) is
\be\label{newMeln}
M^* (A) = \int_{\bf R} ( q_0 (t) - \sin q_0 (t) ) \g (t) +
f ( \om t + A, q_0 (t)) - f ( \om t + A, 0 )
- \partial_q f (\om t + A , 0 ) q_0 (t) \ dt.
\ee
Integrating by parts, since
$ - \ddot{\g} + \g = \partial_q f (\om t + A , 0 ) $, we have
$$
\int_{\bf R} (q_0 (t) - \sin q_0 (t)) \g (t)  \ dt =
\int_{\bf R} ( q_0 (t) - \ddot{q}_0 (t)) \g (t) \ dt =
\int_{\bf R} (- \ddot{\g} (t) + \g (t) ) q_0 (t) dt =
\int_{\bf R}  \partial_q f ( \om t + A , 0 ) q_0 (t) dt,
$$
and we deduce from (\ref{newMeln}) that
$ M^* ( A ) = M ( A ) =
\int_{\bf R} [f( \om t + A, q_0 ( t)) -  f( \om t + A, 0)] \ dt.$
\end{pf}

\begin{theorem}\label{apug}
Assume $(H1)$ and let $ M $ possess a 
proper minimum (or maximum) $ A_0 $, i.e. $ \exists  r >0 $ such that   
$ \inf_{\partial B_r ( A_0 )} \Gamma > \Gamma ( A_0 ) $. 
Then, for $ \mu $ small enough,
the same statement of theorem \ref{thm:main1g} holds where
the diffusion time is $ T_d = O((1/\mu ) \log (1/\mu )).$
\end{theorem}

\begin{remark}
By theorems \ref{thm:invto}-\ref{thm:main1g}
we obtain that, for a priori-stable,
isochronous, degenerate systems considered in \cite{BCV}
$$
{\cal H}_\e = \e \om \cdot I +  \frac{p^2}{2} + \e^d ( \cos q - 1) +
\mu f( \vphi, q ) \quad {\rm with} \quad  1 < d < 2,
$$
for $\mu = \d \e^d $,  $\d $ being a small constant, the diffusion time
is bounded by $ T_d = O (C( \d ) / \e^d )$. This improves
the result of \cite{BCV}, which holds for $ \mu = O ( \e^{d'} )$,
$ d' > d/2 + 3   $, and provides the upper bound on the diffusion time
$ T_d = O( 1 / \e^{C + 2 (\tau + 1)( 2 d' - 1- d/2)})$,
$C $ is a suitable positive constant.
\end{remark}

\section{Splitting of separatrices} \label{justspli}

If the frequency vector $ \om = \om_\e $ contains some
``fast frequencies'' $ \om_i = O( 1/ \e^b )$, $ b > 0 $,  $\e$
being a small parameter, and if the perturbation 
is analytical, the oscillations of
the Melnikov function along some directions turn out to be exponentially
small with respect to $ \e $.
Hence the development  (\ref{eq:eqsvi}) will
provide a valid measure of the splitting only for
$ \mu $ exponentially small with respect to $ \e $.
In order to justify the dominance
of the Poincar\'e-Melnikov function when $ \mu = O(\e^p )$
we need more refined estimates for the error.
The classical way to overcome this difficulty would be to extend
analytically the function $ F_\mu (A, \teta) $ for
complex values of the variables, see \cite{An1}-\cite{DGTJ} and
\cite{Sa}.
However it turns out that the function $F_\mu (A, \teta) $ can 
not be easily analytically extended in a
sufficiently wide complex strip (roughly speaking, the condition
$ q^\mu_{A, \teta } ( Re \ \teta ) = \pi $ appearing
naturally when we try to extend the definition of
$q^\mu_{A, \teta}$ to $\theta \in {\bf C}$ breaks
analyticity). We bypass this problem considering the action functional
evaluated on different ``1-bump pseudo-homoclinic solutions''
$ Q^\mu_{A,\teta} $.
This new ``reduced action functional''
$ {\wtilde F}_\mu (A , \teta) ={\wtilde G}_\mu (A + \om \teta) $
has the advantage to have an analytical extension in 
$ ( A, \teta )$ in a wide complex strip.
Moreover we will show that  the homoclinic functions
$ G_\mu $, $ {\wtilde G}_\mu $ corresponding
to both reductions are the same up to a change of variables
of the torus close to  the identity. This enables
to recover enough information on the
homoclinic function $ G_\mu $ to construct diffusion orbits.
\\[2mm]
\indent
We assume that $ f( \vphi, q ) = ( 1 - \cos q ) f ( \vphi )$,
$ f( \vphi  ) = \sum_{k \in {\bf \footnotesize  Z}^n} f_k \exp({i k \cdot
\vphi})$
and that, there are  $ r_i \geq 0$ such that
\be\label{eq:hypfc}
\forall s \in {\bf N}, \ \exists C_s > 0 \quad {\rm such \ that} \quad
| f_k | \leq \frac{ C_s}{|k|^s} \exp \Big( - \sum_{i=1}^n r_i |k_i|
\Big),  \  \forall k \in {\bf Z}^n.
\ee
Condition (\ref{eq:hypfc}) means that
$ f $ has a $C^{\infty}$ extension defined in
$$
D:=  ( {\bf R} + i [ - r_1, r_1 ]) \times \ldots \times
({\bf R} + i [ -r_n, r_n ])
$$
which is holomorphic w.r.t. the variables for which $r_i > 0$
in $ ( {\bf R} + i I_1) \times \ldots \times
({\bf R} + i I_n)$, where $ I_i=\{0\} $ if $r_i=0$, $I_i
=(-r_i,r_i)$ if $r_i >0$. 
We denote the supremum of $ | f | $ over $ D $ as  
\be\label{eq:supnorm}
||f||:=\sup_{A\in D} |f(A)|.
\ee
It will be used from subsection \ref{secanex}.
\subsection{The change of coordinates}

Define $\psi_0: {\bf R} \to {\bf R}$ by
$ \psi_0 (t) = \cosh^2(t)/ (1+ \cosh t)^3 $ and set
$ \psi_\teta (t) = \psi (t - \teta ).$
Note that
$\int_{{\bf R}} \psi_0 (t) {\dot q}_0 (t) \ dt = \gamma \neq 0$.
Arguing as in lemma \ref{lem:1bump} we can prove

\begin{lemma}\label{lem:2bump}
For $\mu $ small enough (independently of $ \om $),
$ \forall \teta \in {\bf R} $, there exists  a unique function
$ Q_{A, \teta}^\mu (t) : {\bf R} \to {\bf R}$,
and a constant $ \a_{A, \teta}^\mu $ smooth
in $ ( A, \teta, \mu )$, such that

\begin{itemize}
\item \label{eq:secvinr}
$ (i) - \ddot{Q}_{A, \teta}^\mu (t)+ \sin{Q}_{A, \teta}^\mu (t) = \mu \
\sin Q_{A, \teta}^\mu (t) f (\om t + A) +\a_{A, \teta}^\mu \psi_\teta (t);$
\item
$ (ii) \int_{{\bf R}} (Q_{A, \teta}^\mu (t) - q_\teta (t)) \psi_\teta (t) \
dt = 0;$
\item
$ (iii) \max \Big( | Q^\mu_{A, \teta }(t) - q_\teta (t)|,
|{\dot Q}^\mu_{A, \teta }(t) - {\dot q}_\teta (t)| \Big) = O \Big( \mu 
\exp ({- \frac{ |t - \teta |}{2}}) \Big)$;
\item
$ (iv) \max \Big( | \partial_A  Q^\mu_{A, \teta } (t) |,
| \partial_A \dot{Q}^\mu_{A, \teta } (t)|,
|\om . \partial_A  Q^\mu_{A, \teta } (t) |,
|\om . \partial_A \dot{Q}^\mu_{A, \teta } (t)| \Big) =  
O \Big(\mu  \exp({- \frac{ |t - \teta |}{2}}) \Big)$. 
\end{itemize}
Moreover
$ Q^\mu_{A, \teta } (t) = Q^\mu_{A + k 2 \pi, \teta } (t),
\ \forall k \in {\bf Z}^n $ and
$ Q_{A, \teta + \eta }^\mu (t + \eta ) = Q_{A + \om \eta, \teta}^\mu (t)$,
$ \forall \teta, \eta \in {\bf R}$.
\end{lemma}

We define the function
$\wtilde{F}_\mu ( A, \teta):{\bf T}^n \times {\bf R} \to {\bf R} $
as the action functional of Lagrangian \ref{lagraper}
evaluated on the ``1-bump pseudo-homoclinic solutions''
$ Q_{A, \teta }^\mu (t) $ obtained in lemma \ref{lem:2bump}, namely
\be\label{def:homoana}
\wtilde{F}_\mu (A, \teta):=
\int_{\bf R}  {\cal L}_\mu ( Q_{A,\teta}^\mu (t),
\dot Q^\mu_{A,\teta} (t), t) \ dt
\ee
and $ \wtilde{G}_\mu (A):{\bf T}^n  \to {\bf R} $ as
$ \wtilde{G}_\mu ( A ) = \wtilde{F}_\mu ( A, 0 )$.
The following invariance property holds
$\wtilde{F}_\mu (A, \teta + \eta) = \wtilde{F}_\mu
(A + \om \eta, \teta)$, $\forall \teta, \eta \in {\bf R}$; in
particular $\wtilde{F}_\mu (A, \teta ) = \wtilde{G}_\mu (A + \om \teta)$,
$\forall \teta \in {\bf R}$.

\begin{remark}
By lemma \ref{lem:2bump}-(i)-(ii),
if $ \partial_\teta \wtilde{F}_\mu (A, \teta)= 0 $ then
$ Q_{A, \teta }^\mu $ is a true solution of (\ref{pendper}).
More precisely we have
$| \a^\mu_{A, \teta} | \leq C |\partial_\teta \wtilde{F}_\mu (A, \teta) |$,
for a suitable positive contant $C > 0 $.
In addition we could easily prove using lemma \ref{lem:2bump} that
\be\label{eq:estipp}
|\nabla^s \wtilde{G}_\mu( A )  | = O(\mu), \qquad
|\nabla^s \wtilde{F}_\mu (A, \teta) | = O(\mu), \qquad 
 s=1,2.
\ee

\end{remark}

The relation between the two functions $ F_\mu (A, \teta) $ and
$ {\wtilde F}_\mu (A, \teta) $ is given below.
The next theorem is formulated to handle with also 
non-analytical perturbations $ f $. For the analytical 
case see remark \ref{rem:anali}.

\begin{theorem}\label{thm:diffeo}
For $ \mu $ small enough  (independently of $ \om $) 
there exists a Lipschitz continuous function
$ \ov{h}_\mu: {\bf T}^n \times {\bf R} \to {\bf R} $, with
$ \ov{h}_\mu (A, \teta )=O(\mu)$, $| \ov{h}_\mu (A', \teta')-
\ov{h}_\mu (A, \teta) | = O(\sqrt{\mu} (|A'-A|+|\theta'-\theta|) )
$,
$\ov{h}_{\mu}(A,\theta +\eta)= \ov{h}_{\mu} (A+\eta \om, \theta)$,
such that
$ F_\mu ( A, \teta ) = \wtilde{F}_\mu ( A, \teta + \ov{h}_\mu (A, \teta ))$.
In particular, setting $ \ov{g}_\mu (A) = \ov{h}_\mu (A, 0) $,
the homeomorphism $ \psi_\mu : {\bf T}^n \to {\bf T}^n $
given by $\psi_\mu (A) = A +  \ov{g}_\mu (A) \om $ satisfies
$ G_\mu  = \wtilde{G}_\mu  \circ \psi_\mu $.
\end{theorem}

In order to prove theorem \ref{thm:diffeo} we need 
the next two lemmas, proved in the appendix.

\begin{lemma}\label{lem:prdi}
For $ \mu $ small enough (independently of $ \om $)
there exists a smooth function
$ l_\mu (A, \teta) $ with $ l_\mu (A, \teta ) = O ( \mu ),
\nabla l_\mu ( A, \teta ) =  O ( \mu ) $,
$l_{\mu}(A,\theta +\eta)= l_{\mu} (A+\eta \om, \theta)$ such that
$ Q_{A, \teta}^\mu ( \teta + l_\mu (A, \teta)) = \pi.$
\end{lemma}

Define $ V_\mu (A, \teta) = F_\mu (A, \teta + l_\mu (A, \teta) )$.

\begin{lemma}\label{lem:apfuf}
There exists a positive constant $ C_4 $ such that,
for all $(A, \teta) \in {\bf T}^n \times {\bf R} $, there holds 
$$
| \wtilde{F}_\mu ( A, \teta) - V_\mu (A, \teta ) | \leq C_4
| \partial_\teta \wtilde{F}_\mu (A, \teta) |^2.
$$
In particular if $\partial_\teta \wtilde{F}_\mu (A, \teta) = 0 $
then  $ {\wtilde F}_\mu ( A, \teta) = V_\mu (A, \teta ) $.
\end{lemma}

\begin{pfn}{\sc of Theorem} \ref{thm:diffeo}.
By lemma \ref{lem:prdi}, there is a smooth
function $\ov{l}_{\mu}$ such that 
$\ov{l}_{\mu} (A,\theta) = O(\mu)$, 
$\nabla \ov{l}_{\mu} (A,\theta) = O(\mu)$ 
$\ov{l}_{\mu} (A,\theta+\eta)=\ov{l}_{\mu} 
(A + \eta \om, \theta)$ and
$ F_\mu (A,\theta)=V_{\mu}(A,\theta+\ov{l}_{\mu} (A,\theta))$.
 So it is enough to find $ h = h_\mu (A, \teta ) $ such that
\be\label{eq:difeq}
V_\mu (A, \teta ) = \wtilde{F}_\mu ( A, \teta + h ).
\ee
$\ov{h}_{\mu}(A,\theta)$ will be then defined by
\be \label{defovh}
\ov{h}_{\mu}(A,\theta) = \ov{l}_{\mu}(A,\theta) + 
h_{\mu}(A,\theta+ \ov{l}_{\mu}(A,\theta)).
\ee
Note that if $ \partial_\teta \wtilde{F}_\mu (A, \teta ) = 0 $ then, by
lemma
\ref{lem:apfuf}, equation (\ref{eq:difeq}) is solved
by $ h  = 0 $. 
In general we look for $ h $ of the form
$ h = \partial_\teta \wtilde{F}_\mu (A, \teta) g$. Then we can write
\begin{eqnarray}\label{eq:restco}
\wtilde{F}_\mu ( A, \teta + h ) & = &
\wtilde{F}_\mu ( A, \teta ) +  \partial_\teta \wtilde{F}_\mu (A, \teta ) h
+ R_\mu (A, \teta, h ) h^2 \\ \nonumber
& = & \wtilde{F}_\mu ( A, \teta ) +
(\partial_\teta \wtilde{F}_\mu (A, \teta ))^2 g +
R_\mu (A, \teta, \partial_\teta \wtilde{F}_\mu  (A, \teta) g )
(\partial_\teta \wtilde{F}_\mu (A, \teta ))^2 g^2
\end{eqnarray}
where
$$ 
R_\mu (A, \teta, h) = \frac{1}{h^2} \Big[ \wtilde{F}_\mu ( A, \teta + h ) -
\wtilde{F}_\mu ( A, \teta ) -  \partial_\teta
\wtilde{F}_\mu (A, \teta ) h \Big]
$$
is smooth and, by the estimates (\ref{eq:estipp}) 
on the derivatives of $ \wtilde{F}_{\mu}$, it
satisfies $R_\mu (A, \teta, h ) = O(\mu )$,
$\partial_h R_\mu (A, \teta, h ) = O(\mu / |h|)$. 
By (\ref{eq:restco}) equation (\ref{eq:difeq}) is  then equivalent to
$$
\frac{V_\mu (A,\teta) - \wtilde{F}_\mu (A,\teta)}{(\partial_\teta
\wtilde{F}_\mu (A, \teta ))^2 } = g +
R_\mu (A, \teta, \partial_\teta \wtilde{F}_\mu (A, \teta) g ) g^2
$$
We have $R_\mu (A, \teta, \partial_\teta \wtilde{F}_\mu 
(A, \teta) g ) g^2 = O(\mu g^2)$ and 
$\partial_g \Big( R_\mu (A, \teta, \partial_\teta \wtilde{F}_\mu 
(A, \teta) g ) g^2\Big) = O(\mu g)$.
By the contraction mapping theorem, for
$\mu $ small enough, for all $ y \in {\bf R} $ such that
$|y| < 2 C_4$, there exists a unique
solution $ g=\varphi (\mu, A,\theta,y) $ of the equation 
\be \label{eqyy}
y = g +
R_\mu (A, \teta, \partial_\teta \wtilde{F}_\mu (A, \teta) g ) g^2,
\ee
such that $|g| < 3C_4$. Moreover, the function $\varphi$ defined
in this way is smooth.
Setting
\be \label{exph}
h_{\mu} (A, \teta) :=
\varphi \Big(\mu,A,\theta, \frac{V_\mu (A,\teta)
- \wtilde{F}_\mu (A,\teta)}{(\partial_\teta
\wtilde{F}_\mu (A, \teta ))^2 }   \Big) \partial_{\theta} 
\wtilde{F}_{\mu} (A,\theta) \ee 
if $\partial_{\theta} 
\wtilde{F}_{\mu} (A,\theta) \neq 0$ and
$h_{\mu} (A,\theta)=0$ if $\partial_{\theta} 
\wtilde{F}_{\mu} (A,\theta)=0$, 
we get a continuous function $h_{\mu}$ which satisfies 
(\ref{eq:difeq})  and 
$| h_\mu (A, \teta) | \leq
3C_4 | \partial_\teta \wtilde{F}_\mu (A, \teta)| $, which
implies $|h_{\mu}|=O(\mu)$. Moreover 
$h_{\mu}$ is the unique function that enjoys this properties.

By (\ref{exph}) the restriction of $h_{\mu}$ to
$$
U_{\mu}:=\{(A,\theta) \in {\bf T}^n \times {\bf R} \ : \
\partial_\teta \wtilde{F}_\mu (A, \teta) \neq 0  \}
$$
is smooth. Deriving the identity
$ V_\mu (A, \teta)  = {\wtilde F}_\mu (A, \teta +  h_\mu (A, \teta))$
we obtain 
$$
 \partial_\teta {\wtilde F}_\mu (A, \teta +  h_\mu (A, \teta)) 
\nabla h_\mu (A, \teta))=
 \nabla V_\mu (A, \teta) - \nabla 
{\wtilde F}_\mu ( A, \teta +  h_\mu (A, \teta))
$$
for $(A,\theta) \in U_{\mu}$.
Now, from
$$
| \wtilde{F}_\mu ( A, \teta) - V_\mu (A, \teta ) | \leq C_4
| \partial_\teta \wtilde{F}_\mu (A, \teta) |^2, \quad
\partial^2_{\teta \teta} V_\mu  (A, \teta) = O(\mu ), \quad
\partial^2_{\teta \teta} \wtilde{F}_\mu (A, \teta) = O(\mu ),
$$
we can derive  by an elementary argument whose main ingredient 
is Taylor formula that
\be \label{nablaV}
| \nabla \wtilde{F}_\mu ( A, \teta) - \nabla
V_\mu (A, \teta ) | = O( \sqrt{\mu} |\partial_{\theta} \wtilde{F}_\mu ( A,
\teta)|).
\ee
Moreover, by the estimates of $|\partial^2_{\theta \theta}
\wtilde{F}_{\mu}|$ and $|h_{\mu}|$, 
$ (\partial_\teta {\wtilde F}_{\mu})(A, \teta +  h_\mu (A, \teta)) =
(1 + O(\mu) ) \partial_\teta {\wtilde F}_\mu (A, \teta )$.
Hence by (\ref{nablaV}) $|\nabla h_{\mu}|=O(\sqrt{\mu})$ uniformly 
in $U_{\mu}$.
Since $h_{\mu}$ is continuous and $h_{\mu}(A,\theta)=0$ if
$(A,\theta) \notin U_{\mu}$, the Lipschitz continuity of 
$h_{\mu}$ follows.

To complete the proof, we observe that
$ h_\mu (A, \teta +\eta) =  h_\mu (A + \eta \om , \teta)$, which is 
a consequence of  uniqueness. Hence, by (\ref{defovh}) and the properties
of $ \ov{l}_{\mu} $, $ \ov{h}_{\mu} $ satisfies the same.
\end{pfn}

\begin{remark}
Assume that $\wtilde{G}_{\mu}$ satisfies
the ``splitting condition'' \ref{spli} 
(or its generalisation introduced in remark
\ref{respli}), with bounds $\delta$ and
$\alpha$. Then, by theorem \ref{thm:diffeo}, $G_{\mu}$ too
satisfies this condition, with constants $\delta'$, $\alpha'$.
We can take $\delta'=\delta$; moreover, at least if 
$\sqrt{\mu} |\om|$ is small, we can take $\a'=\a/2$.
As a consequence, the results that we shall obtain 
in the next section proving a ``splitting condition'' 
for $ {\wtilde G}_\mu $, may be used to apply the shadowing theorem
\ref{thm:main}.
\end{remark}

\begin{remark}\label{rem:anali}
Assume that $ r_i >0 $ for all $ i $ (i.e. that 
the perturbation $ f $ is analytical). Then we can prove,
using the arguments of the next subsection,
that the homoclinic function $ G_\mu ( \cdot )= F_{\mu}( \cdot , 0 ) $ 
can be extended to a complex analytical function over the interior of $ D $. 
Hence $F_\mu (A, \theta )= G_\mu ( A + \om \theta )$ can be defined
in an open neighbourhood of ${\bf T}^n \times {\bf R}$ in
$({\bf T}^n + i {\bf R}^n ) \times {\bf C}$, so that
the extension is analytical.
One could check that $ l_\mu $ and $ V_\mu $, defined in lemma \ref{lem:prdi}
have analytical extensions too, 
and that the inequality of lemma \ref{lem:apfuf}
still holds in the new set of definition.
Moreover in the next lemma \ref{thm:anaex} it is proved that 
$ \wtilde F_{\mu} $ is analytical w.r.t. $ ( A,\theta ) $.
As a consequence, $(V_\mu (A,\teta) - \wtilde{F}_\mu (A,\teta))/(\partial_\teta
\wtilde{F}_\mu (A, \teta ))^2 $ is real analytical,
and so is the function $ h_{\mu} $ defined in the proof of
theorem \ref{thm:diffeo}. Therefore if $ r_i > 0 $ for all $ i $, then
the homeomorphism $ \psi_{\mu} $ defined in theorem
\ref{thm:diffeo} is a real analytical diffeomorphism. 
\end{remark}

\subsection{Analytical extension} \label{secanex}

The unpertubed homoclinic $ q_0 (t) = 4 {\rm arctg} \ e^t $ 
can be extended to a holomorphic function
over the strip $ S := {\bf R} + i ( -\pi / 2, \pi /2 ) $.
Moreover equation (\ref{pendper}) may be considered also for
complex values of $ q $ and, for $ \mu = 0 $,
$ q_\theta $ is a solution of (\ref{pendper}) for all $ \theta \in S $.
We shall use the notation
$ S_\delta = {\bf R} +
i ( - (\frac{\pi}{ 2} - \delta), \frac{\pi}{2} - \delta) $ ,
for $\delta \in (0, \pi/2)$.
We have
$$
{\dot q}_\teta ( z ) = \frac{2}{\cosh ( z - \teta )},
\quad
{\ddot q}_\teta ( z )= \sin q_\teta (z) = - 2
\frac{\sinh ( z - \teta) }{\cosh^2 ( z - \teta)}, \quad
(1 - \cos q_\teta ( z )) = \frac{2}{\cosh^2 ( z - \teta)}.
$$
Assume that $\theta \in S_{\delta},  Re (\teta) = 0 $.
 The following estimates hold, where $ t \in {\bf R} $.
\begin{eqnarray}
| {\dot q}_\teta ( t )| & \leq & \frac{C}{\min \{(|t|  + \d), 1\}}
\exp ({- |t|}); \label{eq:1}\\
|\sin q_\teta (t)| & \leq & \frac{C}{\min \{(|t|  + \d)^2, 1\}}
\exp({- |t|}); \label{eq:2}\\
|\cos q_\teta (t)| & \leq & \frac{C}{\min \{(|t|  + \d)^2, 1\}};
\label{eq:3} \\
\frac{1}{ | {\dot q}_\teta ( t )| } & \leq & C \exp({|t|})
\min \{ (|t|  + \d), 1\}.  \label{eq:4}
\end{eqnarray}

In what follows
we consider the Banach spaces 
$$
X = \Big\{ w \in C^2 ({\bf R}, {\bf C}) \ \Big|
\ \sup_t  \exp({|t|/2}) (|w(t)| +
|{\dot w}(t)| + |{\ddot w}(t)| ) < + \infty \Big\}
$$
and
$$
\ov{X} = \Big\{ w \in X \ \Big| \ w(0) =0 \Big\},
$$
endowed with norm
$$
||w||_{2,\d} =  \sup_{|t| > 1} \Big( |w(t)| +
|{\dot w}(t)| + |{\ddot w}(t)| \Big) \exp({\frac{|t|}{2}}) +
\sup_{|t| < 1}
\Big( \frac{|w(t)|}{(|t| + \d)^2} +  \frac{| {\dot w}(t)|}{(|t| + \d)} +
|{\ddot w}(t)| \Big).
$$
The next lemma extends lemma \ref{lem:2bump} for
complex values of the variables.
First note that the function $ \psi_0 (t) $ can be extended
to a holomorphic function on ${\bf R} + i( - \pi , \pi) $. 
Recalling the definition for $||f||$ given in (\ref{eq:supnorm}), we have

\begin{lemma}\label{thm:anaex} 
There  exist positive constants $ \eta, C_5 $ such that for all
$ \d \in (0, \pi /2 )$, 
$ \forall 0 <  \mu  \leq (\eta \delta^3)/ || f || $, for all 
$ \om $, for all $ A \in D $,  for all $\teta \in S_\d $ there
exist a unique
$ Q^\mu_{A, \teta}: {\bf R} \to {\bf C}$ and 
a unique $\alpha^\mu_{A, \teta} \in
{\bf C}$ such that
\begin{itemize}
\item
$ Q^\mu_{A, \teta} = q_{\teta + \nu^\mu_{A,\teta}} +
w^\mu_{A, \teta} $, where
 $ \nu^\mu_{A,\teta} \in {\bf C}$, $w\in \ov{X}$ and
$ || w^\mu_{A, \teta}  ||_{2,\d} + | \nu^\mu_{A,\teta} | +
|\alpha^\mu_{A, \teta}|
\leq C_5 \mu ||f|| / \d^2$;
\item
$ - {\ddot Q^\mu_{A, \teta}  }(t) +
\sin{Q^\mu_{A, \teta} }(t) = \mu \ \sin Q^\mu_{A, \teta}(t) \ f ( \om t + A)
+
\a^\mu_{A, \teta} \psi_\teta (t)$;
\item
$\int_{\bf R} (Q^\mu_{A, \teta}(t) - q_{\theta}(t)) \psi_\teta (t)
\ dt = 0$.
\end{itemize}
Moreover
$ Q^\mu_{A, \teta} $ and $\alpha^\mu_{A, \teta} $
depend  analytically on $\theta  $ and on the $A_i$
for which $r_i >0$.
\end{lemma}

\begin{pf}
{\bf 1st step}. Let us  consider  the Banach space 
$$
Y = \Big\{ v \in C ({\bf R}, {\bf C}) \ \Big| \ \sup_t
|v(t)| \exp \Big( {\frac{|t|}{2}} \Big) < + \infty \Big\}
$$
endowed with norm
$ || v ||_{-1,\d} =  \sup_{|t| > 1}  |v(t)| \exp({\frac{|t|}{2}}) +
\sup_{|t| < 1} ( | t | + \d ) | v( t ) |$.
Let $ \teta \in S_\delta $ be given once for all. We may assume without
loss of generality that $Re (\teta ) = 0 $. 

For $ \teta' \in S_{\d/2} $ such that $ | \teta' - \teta | \leq \d/2$ we
introduce the linear operator $ L_{\teta'}: \ov{X} \times {\bf C} \to Y
$
defined by
$$
L_{\teta'} (w,\a) = - {\ddot w} + (\cos q_{\teta'}) w - \a \psi_\teta.
$$
Using that ${\dot q}_{\teta'}$ is a solution of
$ - {\ddot y} + \cos q_{\teta'} y = 0$
we can compute the inverse of $ L_{\teta'} $. It is given by
$ L_{\teta'}^{-1} (g) = (w,\a) $ with
\be
\a  =  -
\frac{\int_{\bf R} g(t) {\dot q}_{\teta'}(t) \ 
dt }{\int_{\bf R} \psi_\teta (t)
{\dot q}_{\teta'} (t) \ dt }  \label{eq:p1},
\ee
\begin{eqnarray}
w(t) & = & {\dot q}_{\teta'} (t) \Big[ \int_0^t -
\frac{1}{{\dot q}_{\teta'}^2 (s)} \Big(   \int_{- \infty}^s
( g( \s ) + \a \psi_\teta ( \s ) ) {\dot q}_{\teta'} ( \s ) \ d \s
\Big) \ ds \Big] \\
& = & {\dot q}_{\teta'} (t) \Big[ \int_0^t
\frac{1}{{\dot q}_{\teta'}^2 (s)} \Big(   \int_s^{+ \infty}
( g(\s) + \a \psi_\teta (\s ) ) {\dot q}_{\teta'} (\s) \ d \s
\Big) \ ds \Big]. \label{eq:p2}
\end{eqnarray}
Note that since $|\teta - \teta'| \leq \d / 2 $,
$ Re ( \teta') \leq \d / 2 $. Therefore estimates
(\ref{eq:1})-(\ref{eq:4}) hold as well (with perhaps
different constants) when $\teta $ is replaced by
$ \teta' $.
We derive from (\ref{eq:p1})-(\ref{eq:p2}) that
\be \label{eq:estinv}
| \a | + || w ||_{2,\d}  \leq \frac{C}{ \d } || g ||_{ - 1, \d }.
\ee

{\bf 2nd Step}. We shall search $ Q $ as $ Q = q_{\teta + \nu} + w $
with $ | \nu | < \d / 2 $, $w \in \overline{X}$. Let $ B $ denote the
open ball
of radius $\d / 2$ in ${\bf C}$ centered at $ 0 $. Let
$J_\mu : B \times \ov{X} \times {\bf C} \to Y  \times {\bf C}$
be defined by
$$
J_\mu ( \nu, w, \a) =
\Big( - \ddot{q}_{\teta + \nu} - \ddot{w}  +
\sin (q_{\teta + \nu} + w) - \mu \sin (q_{\teta + \nu} + w)
\ f ( \vphi ) - \a \psi_\teta,
\int_{\bf R} (q_{\teta + \nu} + w - q_{\teta}) \psi_\teta (t) \Big).
$$
From now we shall use the norms
$ || (\nu, w, \a) ||_2 = |\nu| + ||w||_{2,\delta} + |\a |$
on $ B \times \ov{X} \times {\bf C} $
and
$||(g, \b) ||_{-1} =  ||g||_{-1,\d} + |\b |$ on  $ Y  \times {\bf C}$.
$ J_\mu $ is of class $ C^1 $ and
\begin{eqnarray*}
 D J_\mu (\nu, w, \a) [z,W,a] & = &
\Big( - z \stackrel{\cdots}{q}_{\teta + \nu} - \ddot{W} +
\cos ( q_{\teta + \nu} + w) (z
{\dot q}_{\teta + \nu} + W ) \\
& - & \mu \cos (q_{\teta + \nu} + w)  (z {\dot q}_{\teta + \nu} + W)
\ f ( \vphi ) - a \psi_\teta,
\int_{{\bf R}} (z {\dot q}_{\teta + \nu} + W ) \psi_\teta (t) \Big).
\end{eqnarray*}
We shall prove that, provided $||(\nu, w,\a)||_2 / \d $
and $\mu ||f|| / \d^3 $ are small enough $ DJ_\mu (\nu, w,\a) $
is invertible. We first consider the case when
$ w = 0 $ and $ \mu = 0 $.
Let $ T_\nu = D J_0 (\nu, 0,\a) $ (independent of $\a$). Observing
that $ - \stackrel{\cdots}{q}_{\teta + \nu} +
\cos ( q_{\teta + \nu} )
{\dot q}_{\teta + \nu} =0$, we obtain
$$
T_\nu [z, W, a] = \Big( - \ddot{W} +
\cos{q_{\teta + \nu} } W - a \psi_\teta,
\int_{\bf R} (z {\dot q}_{\teta + \nu} + W ) \psi_\teta (t) \Big).
$$
Using the first step we derive that $ T_\nu $ is invertible
and, for a suitable positive constant $ C $
\be \label{invTnu}
|| T_\nu^{-1} ( g, \b ) ||_2 \leq \frac{C}{\d} ||(g, \b)||_{-1}.
\ee
Now we estimate $ || (D J_\mu (\nu, w, \a) -  T_\nu)[z, W,a]||_{-1}.$
We have
\begin{eqnarray*}
(D J_\mu (\nu, w, \a) -  T_\nu) [z, W,a] & = &
\Big( (\cos (q_{\teta + \nu } + w ) - (\cos q_{\teta + \nu} ) )
(z {\dot q}_{\teta + \nu} + W) \\
& - & \mu \cos (q_{\teta + \nu} + w )
(z {\dot q}_{\teta + \nu} + W) f (\vphi), 0 \Big).
\end{eqnarray*}
We easily get
\begin{eqnarray*}
|| (D J_\mu (\nu, w, \a) -  T_\nu)[z, W,a]||_{-1} & \leq &
C ||w||_{2,\d}  \Big( ||W||_{2,\d} + | z | \Big)
+ \frac{\mu||f||}{\d^2} |z|  +
| \mu | ||f|| ||W||_{2,\d} \\
& \leq & C \Big( ||w||_{2,\d} + \frac{\mu ||f||}{\d^2} \Big)
||(z,W,a)||_2.
\end{eqnarray*}
As a consequence, recalling (\ref{invTnu}), if  
$ \mu ||f|| / \d^3 \leq K_0 $ and 
$ || w||_{2,\d} / \d \leq K_0 $, 
for $ K_0 $ small enough, then $ D J_\mu (\mu, w,\a) $
is invertible and
$$
|| (D J_\mu (\nu, w,\a))^{-1} || \leq \frac{K_1}{ \d }
$$
for a suitable positive constant $ K_1 $.

{\bf 3rd Step.}
We now prove the existence of a constant $K_2$
such that $(0,0,0)$ is the unique solution
of the equation $J_0(\nu,w,\a)=0$ in $B(K_2 \delta)$,
ball centered at the origin and of radius $K_2 \delta$
for the norm $|| \  ||_2$. First we observe that, 
since $\ddot{q}_{\theta +\nu} = \sin(q_{\theta + \nu})$, there holds
$$
J_0 (\nu,w,\a) = T_{\nu}[\nu,w,\a] + \Big( \sin(q_{\theta +\nu} +w) -
\sin(q_{\theta +\nu} ) - \cos(q_{\theta +\nu} ) w , 
\int_{\bf R} (q_{\theta +\nu} -q_{\theta}-\nu \dot{q}_{\theta +\nu})
\psi_{\theta} \Big).
$$
Moreover, by the analyticity of $ q_0,\dot{q}_0 , \psi_0 $ over $ S $,
there holds
$$
\int_{\bf R} (q_{\theta +\nu}(t) -q_{\theta}(t)-\nu \dot{q}_{\theta +\nu}(t))
\psi_{\theta}(t) \ dt  = 
\int_{\bf R} (q_{\nu}(t) -q_{0}(t) - \nu \dot{q}_{\nu} (t))
\psi_{0} (t) \ dt,
$$
hence there is a constant $C'$ such that
$$
\Big| \Big| \Big( \sin(q_{\theta +\nu} +w) -
\sin(q_{\theta +\nu} ) - \cos(q_{\theta +\nu} ) w , 
\int_{\bf R} (q_{\theta +\nu} -q_{\theta}-\nu \dot{q}_{\theta +\nu})
\psi_{\theta} \Big)  \Big| \Big|_{-1} \leq C' (||w||_{2,\delta}^2 + |\nu|^2).
$$
So, if $ J_0(\nu,w,\a)=0 $ then, by (\ref{invTnu})
\begin{eqnarray*} 
||(\nu,w,\a)||_2& =& \Big| \Big| - T_{\nu}^{-1} \Big( 
\sin(q_{\theta +\nu} +w) -
\sin(q_{\theta +\nu} ) - \cos(q_{\theta +\nu} ) w , 
\int_{\bf R} (q_{\theta +\nu} -q_{\theta}-\nu \dot{q}_{\theta +\nu})
\psi_{\theta} \Big) \Big| \Big|_2 \\
&\leq& \frac{CC'}{\delta} ||(\nu,w,\a)||_2^2. 
\end{eqnarray*}
Let $K_2 < 1/(CC')$. By the latter inequality,
if $J_0 (\nu,w,\a)=0$ and $||(\nu,w,\a)||_2 \leq K_2 \d$, then
$\nu=0$, $w=0$, $\alpha =0$.  
\\[2mm]
\indent {\bf 4th step.} By the previous steps
we know that there exist positive constants 
$ K_0 $, $ K_1 $ and $ K_2 $ such that
\begin{itemize}
\item $(i)$ ( $J_0 (\nu,w,\a)=0$ and $||(\nu,w,\a)||_2 \leq K_2 \d$) $\iff$
$\nu=w=\a=0$;
\item $(ii)$ If $| \nu | <  \d / 2 $,  $|| w ||_{2,\d} \leq K_0 \d $,
$\mu ||f|| \leq K_0 \d^3 $  then
$D J_\mu (\nu, w,\a) $ is invertible and
$|| (D J_\mu (\nu, w,\a))^{-1} || \leq K_1 / \d $.
\end{itemize}
Moreover there exists a constant $ K_3 > 0 $ such that
\begin{itemize}
\item $(iii)$
$ || \partial_{\mu} J_\mu (\nu, w,\a)||
_{-1} =
\Big | \Big| (\sin (q_{\teta + \nu} + w) f( \vphi ),0 ) \Big| \Big|_{-1}
\leq ||f|| K_3/ \d$.
\end{itemize}
We say that $(i),(ii),(iii)$ imply that
there is $ \eta $ such that,
for all $0 < \mu  < \eta \delta^3 / ||f||$,
the equation $J_{\mu} (\nu,w,\a)=0$ has a unique solution
such that $||(\nu,w,\a)||_{2}< K_2 \d/2$. In addition
$||(\nu,w,\a)||_2 = O(\mu ||f|| / \d^2)$.
To prove  existence, we can proceed as follows.
Let $\cal S$ denote the set of all $\mu \in [0,K_0 \d^3 / ||f||]$
such that there exists a $C^1$ function
$ S_{\mu} : [0,\mu] \to \{ (\nu,w,\a) \ : \  
||(\nu,w,\a)||_{2} < K_2 \d / 2 \}$
such that $S_{\mu} (0)=0$, $ J_t (S_{\mu}(t))=0 $ for all $t\in [0,\mu]$.
$\cal S$ is a bounded interval. Let us call $\ov{\mu}$
its supremum. By $(ii)$ and  the Implicit Function Theorem,
$\ov{\mu} >0$. In addition, for $\mu \in {\cal S}$, there is a unique
function $ S_{\mu} $  with 
the required properties. As a consequence, for $ 0 < \mu < \mu' $,
$ S_{\mu}=S_{\mu'|[0,\mu]}$ and we can define a $ C^1 $ function
$ S: [0,\ov{\mu}) \to \{ (\nu,w,\a) \ : \  ||(\nu,w,\a)||_2 < 
K_2 \d / 2 \}$ such that $S(t)=S_{\mu}(t)$ for all $\mu \in (0,\ov{\mu})$.
By $(ii)$ and $(iii)$, we can write,  for all
$t\in (0, \ov{\mu})$,
$$ 
||S'(t)||_2 = \Big| \Big| \Big[ DJ_t (S(t)) \Big]^{-1}.
\Big( \frac{\partial J_t}{\partial t} (S(t))\Big) \Big| \Big|_2
\leq \frac{K_1 K_3 ||f||}{\delta^2}.
$$
Hence
\be \label{Sbound}
||S(t)||_2 \leq \frac{K_1K_3}{\d^2}||f|| |t|.
\ee
Now, since $S'(t)$ is bounded, $S(t)$ converges to some
$\ov{S}$ as $t\to \ov{\mu}$. Either $\ov{\mu} = K_0 \delta^3 / ||f||$
or $||\ov{S}||_2 = K_2 \d /2$  (If not, by the 
Implicit Function Theorem, we could extend the solution $S$
to an interval $[0, \ov{\mu }+ \xi)$, $ \xi>0 $, contradicting
the definition of $\ov{\mu}$).
In the latter case, by (\ref{Sbound}),
$$
\ov{S}=\frac{K_2 \d}{2} \leq \frac{K_1 K_3}{\d^2} \ov{\mu}||f||.
$$
So the existence assertion holds for $ 0 < \mu < \eta \d^3 / ||f||$, 
where $ \eta = {\rm min}(K_0 , K_2/(2K_1 K_3))$.

In order to prove uniqueness, we assume
that there are $b_1,b_2$ such that
$||b_i||_2 < K_2 \delta /2$, $J_{\mu} (b_i)=0$.
Then, by the same argument as previously,
we can prove the existence of two functions of class $C^1$
$ S_1, S_2 \ : \  [0,\mu] \to
\{ b \ : \ ||b||_2 < K_2 \d \} $
such that $ S_i(\mu)=b_i $, $ J_t(S_i(t))=0 $.
Moreover, by $(ii)$ and  the Implicit Function Theorem,
 $S_1(\mu)\neq S_2(\mu)$ implies that
$S_1(t)\neq S_2(t)$ for all $t\in [0,\mu]$, which contradicts
$(i)$, proving uniqueness.
 
The bound of
$||w_{A,\theta}^{\mu}||_{2,\d}+|\nu_{A,\theta}^{\mu}|
+|\alpha_{A,\theta}^{\mu}|$
given in the statement is a direct consequence of (\ref{Sbound}).

To complete the proof, we point out that 
$J_{\mu}$ is analytical on $(A,\theta)$.
Therefore, as a consequence of the Implicit Function Theorem (see for example
\cite{An1}),
$$
Q_{A,\teta}^\mu = q_{\teta + \nu^\mu (A,\teta)} + w^\mu(A,\teta)
$$ 
depends analytically on $  \teta $ and on $A_i$ if $r_i>0$.
\end{pf}

We now consider the analytical extension of the
function $ \wtilde{F}_\mu (A, \teta ) $
for $ ( A, \teta ) \in D \times  S_\d $
$$
{\wtilde F}_\mu (A, \teta )  = \int_{\bf R}
\frac{ ({\dot Q}^\mu_{A,\teta})^2 (t)}{ 2 } +
( 1- \cos  Q^\mu_{A,\teta} (t)) +
\mu (\cos  Q^\mu_{A,\teta}(t) - 1 ) f(\om t + A) \ dt.
$$
Let consider also the analytical extension for
$ (A, \teta ) \in D \times S_\d $ of the Melnikov function
$$
M(A, \teta) =  \int_{\bf R} ( 1 - \cos q_\teta (t)))
f ( \om t + A  )  \  dt.
$$
We have $ \Gamma (A + \om \teta)= M (A, \teta) $. We now prove
\begin{lemma}
For $ \mu || f || \delta^{-3} $ small enough,
for all $ ( A, \teta) \in D \times S_\d $, we have
\be\label{eq:stimana}
{\wtilde F}_\mu (A, \teta ) = const +
\mu M (A , \theta ) + O \Big( \frac{\mu^2 ||f||^2 }{\d^4} \Big).
\ee
\end{lemma}

\begin{pf}
We have $ Q^\mu_{A,\teta} = q_{\teta + \nu^\mu_{A,\teta}} + 
w^\mu_{A,\teta}$ and we set for brevity 
$ Q^\mu_{A,\teta}= q_{\teta + \nu } + w $. 
\begin{eqnarray*}
{\wtilde F}_\mu ( A, \teta ) & = &
\int_{\bf R} \frac{ { \dot q_{\teta + \nu} + \dot w}^2}{2}
+ ( 1 - \cos (q_{\teta + \nu} + w) ) +
\mu ( \cos  ( q_{\teta + \nu } + w) - 1 ) f(\om t + A) \ dt \\
& = &
const + \int_{\bf R} - {\ddot q}_{\teta + \nu} w + \frac{1}{2} \dot{w}^2 +
( \cos q_{\teta + \nu} - \cos (q_{\teta + \nu} + w )) \\
& + & \mu (1- \cos  (q_{\teta + \nu} ) f(\om t + A) +
\mu ( \cos q_{\teta + \nu} - \cos (q_{\teta + \nu} + w)) f(\om t + A) \\
& = & const + \mu M ( \teta + \nu, A) + \int_{\bf R}
\frac{1}{2} \dot{w}^2 + \Big( \cos q_{\teta + \nu} -
\cos (q_{\teta + \nu} + w) - \sin q_{\teta + \nu} w \Big) \\
& + & \mu \int_{\bf R}
\Big( \cos q_{\teta + \nu} - 
\cos (q_{\teta + \nu} + w) \Big) f ( \om t + A ).
\end{eqnarray*}
By the estimate $ || w||_{2,\d} \leq C \mu ||f|| / \d^2 $, it 
follows easily
$$
{\wtilde F}_\mu (\teta, A) =
Const + \mu M(\teta + \nu, A) + O \Big( 
\frac{\mu^2 ||f||^2 }{\d^4} \Big).
$$
For example we can get that
$ \int_{\bf R} \cos q_{\teta + \nu} - \cos (q_{\teta + \nu} + w)
- (\sin  q_{\teta + \nu}) w
= O( \mu^2 ||f||^2 / \d^4 )$ by writing
$ \cos q_{\teta + \nu} - \cos (q_{\teta + \nu} + w)
- (\sin  q_{\teta + \nu}) w = w^2 \int_0^1
- (1 - s ) \cos (q_{\teta + \nu} + sw ) \ d s $
and using (\ref{eq:2})-(\ref{eq:3}) togheter with
$|| w ||_{2, \d} \leq \mu ||f|| / \d^2 $.
Moreover
$$
| M ( \teta + \nu , A ) - M ( \teta , A)| = O \Big( \frac{|\nu|}{\d^2}\Big)
=
O \Big( \frac{\mu ||f|| }{\d^4} \Big),
$$
which completes the proof of the lemma.
\end{pf}

The Fourier coefficients of the Melnikov function
$ \Gamma (A) = \sum_k \Gamma_k \exp({i k A })$ are explicitely given by
\be\label{eq:melcoe}
\Gamma_k = f_k \frac{2 \pi (k \cdot \om )}{\sinh ((k \cdot \om)
\frac{\pi}{2}) }.
\ee
By estimate (\ref{eq:stimana}), since ${\wtilde F}_\mu (A, \teta ) =
{\wtilde G}_\mu (A + \om \teta ) $ and
$ M (A, \teta ) =  \Gamma (A + \om \teta)$, via a standard lemma
on Fourier coefficients of analytical functions
(lemma 3 in \cite{DGTJ}), we obtain the following result
(compare with theorem 3.4.5 in \cite{LMS}).

\begin{theorem}\label{thm:tfc}
There exists a positive constant $ C_6 $ such that, for 
$ \mu  ||f|| \d^{-3} $ small enough, then $ \forall k \neq 0,
k \in  {\bf Z}^n$, for all $ \d \in (0, \frac{\pi}{2}) $, for all 
$\om $,
\be\label{eq:devecoe}
| \wtilde{G}_k - \mu \Gamma_k | \leq
\frac{C_6 \mu^2 || f ||^2 }{ \d^4}
\exp \Big( - \sum_{i=1}^n  r_i | k_i | \Big)
\exp \Big( - | k \cdot \om | \Big( \frac{\pi}{2} - \d \Big)  \Big).
\ee
\end{theorem}

\section{Three time scales}

We consider in this section three time  scales systems as
(see \cite{GGM} and \cite{PV})
$$
{\cal H} = \frac{I_1}{\sqrt{\e}}+ \e^a \beta \cdot I_2  + \frac{p^2}{2}
+ (\cos q - 1 ) + \mu (\cos q - 1 ) f(\varphi_1,\varphi_2 ), \quad 
\e > 0
$$
with $ n\geq 2 $,  $\vphi_1 \in {\bf T}^1$,  $\varphi_2 \in {\bf T}^{n-1}$,
$ I_1 \in {\bf R}^1 $,  $ I_2 \in {\bf R}^{n-1} $, $ \b \in {\bf
R}^{n-1} $ and $ \e $ is a positive small parameter.
The frequency vector is $ \om=(1/\sqrt{\e}, \e^a \beta) $, where
$ \beta=(\beta_2, \ldots, \beta_n) \in {\bf R}^{n-1}$ is given.
\\[1mm]
We assume through this section that 
$ \mu || f|| \e^{-3/2}$ and $ \e $ are small.
\\[1mm]
Given $ \kappa_2 = (k_2, \ldots , k_n) \in {\bf Z}^{n-1}$ , 
we shall use the notation $ \kappa_2^+ := (|k_2|, \ldots, |k_n|)$. 
Moreover we shall use the 
abbreviation $ \rho_2 := (r_2, \ldots, r_n) $,    
so that $ \kappa_2^+ \cdot \rho_2 := \sum_{i=2}^n r_i |k_i| $. 
We recall that $ r_1, \ldots, r_n $ 
are defined in formula (\ref{eq:hypfc}).

Writing
$$
f ( \vphi_1, \vphi_2) = \sum_{(k_1,\kappa_2) \in {\bf Z} \times {\bf Z}^{n-1}}
f_{k_1,\kappa_2} \exp({i ( k_1 \varphi_1 + \kappa_2 \cdot \varphi_2 )}),
$$
we assume  that  $ f $ is analytical w.r.t $ \vphi_2 $.
More precisely, $ r_1=0 $ and for $ i\geq 2 $, $ r_i>0 $. 
If  $ a = 0 $, we impose in addition that $ r_i > |\beta_i| \pi /2  $  for $i\geq 2$.
\\[1mm]
\indent
We shall use  (\ref{eq:devecoe}) in order to give an expansion for the
``homoclinic function''
$$
\wtilde{G}_\mu (A) = \sum_{ (k_1,\kappa_2) \in {\bf Z} \times {\bf Z}^{n-1} }
\wtilde{G}_{k_1,\kappa_2} \exp({i ( k_1 A_1 + \kappa_2 \cdot  A_2 ) }) =
\sum_{k_1 \in {\bf Z}} \wtilde{g}_{k_1} (A_2 ) \exp({i  k_1 A_1}).
$$
We start with

\begin{lemma} \label{lem:3tsbound}
There exists a positive constant $ C_7 $ such that,
for  $\mu ||f||\e^{- 3/2}$ small enough, 
\be\label{eq:restk2}
\sum_{\kappa_2 \in {\bf Z}^{n-1}, | k_1 | \geq 2}
| {\wtilde G}_{k_1,\kappa_2} | \leq C_7
\frac{\mu ||f|| }{\sqrt{\e}} \exp({- \frac{ \pi }{ \sqrt{\e }}}).
\ee
\end{lemma}

\begin{pf}
Choosing  $ \d = \sqrt{\e} $, we get from (\ref{eq:melcoe}) and
(\ref{eq:devecoe})
\begin{eqnarray*}
| \wtilde{G}_{k_1,\kappa_2} | &\leq & \mu |\Gamma_{k_1,\kappa_2}|
+ |\wtilde{G}_{k_1,\kappa_2} - \mu \Gamma_{k_1,\kappa_2}| \\ &
\leq& 
 C \mu ||f|| e^{-\kappa_2^+ \cdot \rho_2 } 
\Big( \Big| \frac{k_1}{\sqrt{\e}}+ \kappa_2 \cdot \beta \e^a \Big| 
+1 \Big) 
e^{-|\frac{k_1}{\sqrt{\e}}+ \kappa_2 \cdot \beta \e^a|\pi /2} \\
&+& C \frac{\mu^2}{\e^2} ||f||^2 e^{-\kappa_2^+ \cdot \rho_2 }
e^{-|\frac{k_1}{\sqrt{\e}}+ \kappa_2\cdot \beta \e^a|(\pi /2-\sqrt{\e})}
\\
&\leq&  C
\frac{\mu ||f|| }{\sqrt{\e}} (|k_1|+|\kappa_2|)
e^{-\kappa_2^+ \cdot \rho_2   + |\kappa_2 \cdot \beta| \e^a \pi /2}
e^{-\frac{|k_1|}{\sqrt{\e} } \pi /2 } \\
&+& C \frac{\mu^2}{\e^2} ||f||^2
e^{-\kappa_2^+ \cdot \rho_2 + |\kappa_2 \cdot \beta|  \e^a (\pi /2 - \sqrt{ \e })}
e^{-|k_1| (\frac{\pi}{2 \sqrt{\e}}-1)} \\
& \leq & C \frac{\mu ||f|| }{\sqrt{\e}} - \pi
 (|k_1|+|\kappa_2|) \exp (- \sum_{j=2}^n |k_j|(r_j - |\beta_j|
\e^a  \pi/2 )) \exp (-|k_1| (\frac{\pi}{2 \sqrt{\e}}-1)). 
\end{eqnarray*}
We have used in the last line that $\mu ||f|| / \e^{3/2} = O(1)$.
Now  $ r_j - |\beta_j| \e^a  \pi / 2  > 0 $ for $\e$ small enough
both if $ a=0$ or if $ a > 0$.
Summing in $ | k_1 | > 2 $ and in $ \kappa_2 \in {\bf Z}^{n-1} $ 
we obtain (\ref{eq:restk2}).
\end{pf}

Let
$$
\Gamma (\e, A) = \sum_{(k_1,\kappa_2) \in {\bf Z} \times {\bf Z}^{n-1} }
\Gamma_{k_1,\kappa_2} \exp({i ( k_1 A_1 + \kappa_2 \cdot  A_2 ) }) =
\sum_{k_1 \in {\bf Z}} \Gamma_{k_1} (\e, A_2) \exp({i  k_1 A_1}).
$$

\begin{lemma}
We have
$$
\wtilde{g}_0 (A_2 ) = \mu \Gamma_0 (\e,  A _2) + O  ( \mu^2 ||f||^2 ).
$$
\end{lemma}

\begin{pf}
A summation over $\kappa_2$ in  estimate (\ref{eq:devecoe}) (where
we chose  $\delta = \pi/2$ and $k_1=0$) yields immediately the estimate.
\end{pf}

\begin{lemma}\label{lem:tre}
We have
$$
\wtilde{g}_{\pm 1} (A_2 ) = \mu \Gamma_ {\pm 1} (\e, A _2) +
O  \Big( \frac{\mu^2 ||f||^2}{\e^2} 
\exp({- \frac{\pi}{2 \sqrt{\e}}}) \Big).
$$
\end{lemma}
\begin{pf}
By (\ref{eq:devecoe}) (where we chose $\delta =\sqrt{\e}$ and  $ k_1 =
\pm 1$), we can obtain as in the proof of lemma \ref{lem:3tsbound}
\begin{eqnarray*}
|\wtilde{g}_{\pm 1}(A_2)-\mu \Gamma_{\pm 1}(\e,  A_2)| & \leq &
C\frac{\mu^2}{\e^2} ||f||^2 \sum_{\kappa_2 \in {\bf Z}^{n-1}}
\exp ({-\sum_j^n |k_j|(r_j - |\beta_j| \e^a \pi /2 ))}
\exp ({- (\frac{\pi}{2 \sqrt{\e}}-1)}.
\\
&\leq & C\frac{\mu^2}{\e^2} ||f||^2
e^{-\frac{\pi}{2 \sqrt{\e}}}.
\end{eqnarray*}
\end{pf}

Since $ \Gamma (A) $ and $ G_\mu (A) $ are real functions
we have that $ \wtilde{g}_{-1} ( A_2 )= 
\ov{\wtilde{g}}_1 ( A_2 ) $ and
$ \Gamma_{-1} ( A_2 ) = \ov{\Gamma}_1 ( A_2 ) $, 
where $ \ov{z} $ denotes the complex conjugate of
the complex number $ z $. We deduce
from  the previous three lemmas the following result.

\begin{theorem}\label{thm:tts}
For $\mu ||f|| \e^{- 3/2} $ small there holds
\begin{eqnarray*}
\wtilde{G}_\mu (A_1, A_2) & = & Const +
\Big( \mu \Gamma_0 (\e,  A _2) + R_0 (\e, \mu,  A_2) \Big) +
2 {\rm Re} \
\Big[ \Big( \mu \Gamma_1 (\e, A _2)
+ R_1 (\e, \mu, A_2) \Big) e^{i A_1} \Big]\\
& + &
O (\mu  \e^{-1/2} ||f|| \exp({ - \frac{\pi}{\sqrt{\e}} })  )
\end{eqnarray*}
where
$$
R_0 (\e, \mu,  A_2) = O \Big( \mu^2 ||f||^2 \Big) \quad {\rm and}
\quad R_1 (\e, \mu, A_2) = O \Big( \frac{\mu^2 ||f||^2}{\e^2}
\exp({- \frac{\pi}{2 \sqrt{\e}}}) \Big).
$$
\end{theorem}

\begin{remark}\label{rem:tts}
(i) This improves the results in \cite{PV} which require
$ \mu = \e^p $ with $ p > 2 + a $.

(ii) Theorem \ref{thm:tts} certainly holds in any dimension,
while the results of \cite{GGM}, which hold for more general systems,
are proved for 2 rotators only.

(iii) Theorem \ref{thm:tts} is not in contradiction
with the counterexample given in \cite{GGM3}.

(iv) In order to prove a splitting condition using theorem \ref{thm:tts}
it is necessary, according with \cite{GGM} and \cite{PV},
that $ \exists m, l \in {\bf Z}^{n-1}$ such that
$ f_{0,l}, f_{1, m} \neq 0 $. Otherwise, recalling 
(\ref{eq:melcoe}), it results that
$ \Gamma_0 (\e, A_2) = \sum_{\kappa_2 \in {\bf Z}^{n-1}} 
\Gamma_{0, \kappa_2} \exp^{ i \kappa_2 \cdot A_2} = 0$ and also
$ \Gamma_1 (\e, A_2) = \sum_{\kappa_2 \in {\bf Z}^{n-1}}
 \Gamma_{1, \kappa_2} \exp^{ i \kappa_2 \cdot A_2} = 0 $.
\end{remark}

Theorem \ref{thm:tts} enables us to provide
conditions  implying the existence of diffusion orbits.
For instance we obtain the following result.

\begin{lemma} \label{lem:3ts}
Assume that there are $ \ov{A}_2 \in {\bf R}^{n-1}$ 
and $ d_0,c_0 > 0 $ such that, for all small $ \e >0 $,
$$
\begin{array}{rl}
(i) & \quad  |\Gamma_1 ( \e, A_2 )| > (c_0/\sqrt{\e}) e^{-\pi /
( 2 \sqrt{\e})}, 
\quad \forall A_2 \in {\bf R}^{n-1} \ \ {\rm such \  that} \ \
|A_2 - \overline{A}_2| \leq  d_0; \\
(ii) & \quad \Gamma_0 (\e, A_2 ) > \Gamma_0 ( \e, \ov{A}_2) + c_0,
\quad \forall A_2 \in {\bf R}^{n-1}, \ \ {\rm such \  that} \ \
|A_2 - \overline{A}_2| =  d_0.  
\end{array}
$$
Then there is $c_1 >0$ such that, 
for $\mu ||f|| \e^{-3/2}$ small enough, condition \ref{spli}
is satisfied by $ {\wtilde G}_\mu $, with
$ \alpha = c_1 e^{-\pi / (2 \sqrt{\e})} / \sqrt{\e}$ and  $ \delta =
c_0\mu/(2\sqrt{\e}) e^{-\pi / (2 \sqrt{\e})}$.
\end{lemma}

\begin{pf}
First we can derive from (\ref{eq:melcoe}) and (\ref{eq:devecoe})
in the same way as in the proof of lemmas 
\ref{lem:3tsbound} and \ref{lem:tre} that
\be \label{g1b}
\begin{array}{rlc}
|\wtilde{g}_1 (A_2)| + |\nabla \wtilde{g}_1 (A_2) | & \leq &
\dps \sum_{\kappa_2 \in {\bf Z}^{n-1}} (1+|\kappa_2|)
|\wtilde{G}_{1,\kappa_2}| \\  \\
&\leq& C \dps \frac{\mu ||f||}{\sqrt{\e}} e^{- \pi /(2\sqrt{\e})}.
\end{array}
\ee
By the bounds of $ R_0 $  and $ R_1 $ of theorem \ref{thm:tts},
for $\e$ and  $\mu ||f|| \e^{-3/2}$ small enough, we have
$$
\begin{array}{rl}
(i) & \quad  |\wtilde{g}_1 (A_2)|= |(\mu \Gamma_1+ R_1)  ( A_2 )|>
(\mu c_0/(2\sqrt{\e})) e^{-\pi /( 2 \sqrt{\e})} 
\quad \forall A_2 \in B_{d_0} \\
(ii) & \quad  |\wtilde{g}_0 (A_2)|=
(\mu \Gamma_0 +R_0)({A}_2 ) > (\mu \Gamma_0 + R_0) (\ov{A}_2) +c_0/2
\quad \forall A_2 \in \partial B_{d_0},
\end{array}
$$
where $B_{d_0}$ is the open ball centered at $\ov{A}_2$ 
of radius $d_0$.

So we can write $\wtilde{g}_1 (A_2) = |\wtilde{g}_1 (A_2)| e^{i
\phi(A_2)}$, where $\phi$ is a smooth function defined in
$B_{d_0}$ ( (\ref{g1b}) and the previous lower bound
of $|\wtilde{g}_1(A_2)|$ provide a bound of $\nabla \phi (A_2)$).

For $A_2 \in B_{d_0}$, by theorem \ref{thm:tts} we have
$$
\wtilde{G}_{\mu} (A_1,A_2)= Const + (\mu \Gamma_0+ R_0) (\e,\mu,A_2) +
2|(\mu \Gamma_1+ R_1) (\e,\mu,A_2)| \cos(A_1+\phi(A_2))
+ O(\mu \e^{-1/2}||f|| \exp(-\frac{\pi}{\sqrt{\e}})).
$$
Let
$$
U=\{(A_1,A_2) \in {\bf R} \times {\bf R}^{n-1} \ : \  A_2 \in B_{d_0}, \
|A_1 + \phi(A_2) - \pi | < \frac{\pi}{2}  \}.
$$
If $A_2 \in \partial B_{d_0}$ then
$$
\wtilde{G}_{\mu} (A_1,A_2) - \wtilde{G}_{\mu}
(\pi-\phi(\ov{A}_2),\ov{A}_2)
\geq \frac{c_0}{2}+ O(\mu \e^{-1/2}||f|| \exp(-\frac{\pi}{2\sqrt{\e}})).
$$
If $|A_1+\phi(A_2)-\pi| =\pi/2$ then
\begin{eqnarray*}
\wtilde{G}_{\mu} (A_1,A_2) - \wtilde{G}_{\mu}
(\pi - \phi( A_2 ), A_2) & = & 2|(\mu \Gamma_1+ R_1)(A_2)| + 
O \Big( \mu \e^{-1/2}||f|| \exp(-\frac{\pi}{\sqrt{\e}}) \Big) \\
& \geq &  \frac{\mu c_0}{\sqrt{\e}} e^{-\pi /( 2 \sqrt{\e})}+
O(\mu \e^{-1/2}||f|| \exp(-\frac{\pi}{\sqrt{\e}})).
\end{eqnarray*}
Hence, for $\e$ and $\mu || f ||\e^{-3/2}$ small enough,
$$
\inf_{\partial U} \wtilde{G}_{\mu} >
\inf_{U} \wtilde{G}_{\mu} + 
\frac{c_0 \mu e^{-\pi / 2 \sqrt{\e}}}{2\sqrt{\e}}.
$$
Using that $|\nabla \wtilde{G}_{\mu}|  = O(\mu)$,
we can easily derive that condition \ref{spli}
(not with a ball $B_{\rho}$ but the bounded open set $U$,
according to remark \ref{respli}) is satisfied
with 
$\delta =  (c_0/2) \mu e^{-\pi / 2 \sqrt{\e}} / \sqrt{\e}$,
$\alpha = c_1 e^{-\pi / 2 \sqrt{\e}} / \sqrt{\e}$ for some 
positive constant $c_1$.
\end{pf}

The condition given in the previous lemma is not easily
handable. We now want to provide simpler conditions,
involving properties of the perturbation $ f $.
For $ A = ( A_1, A_2 ) \in {\bf T}^1 \times {\bf T}^{n-1} $, let
$$
f ( A_1, A_2 ) = 
\sum_{ (k_1, \kappa_2) \in {\bf Z} \times {\bf Z}^{n-1} }
f_{ k_1, \kappa_2 } \exp({i ( k_1 A_1 + \kappa_2 \cdot  A_2 ) }) =
\sum_{k_1 \in {\bf Z}} f_{k_1} (A_2 ) \exp({i  k_1 A_1}).
$$
If $ f $ is analytical in some domain
then also $ f_{k_1}( \cdot ) $ can be analitically extended
is the same domain, as
$ f_{k_1} ( s ) = (1 / 2 \pi ) 
\int_0^{2 \pi} f( \sigma, s) e^{-ik_1\sigma} \ d\sigma.$

\begin{theorem} \label{thm:3tsm}
Assume that $ f $ satisfies one of the following conditions:
\\[1mm]
($i$) $ a>0 $, $f_0 (A_2)$ admits a strict local minimum at
the point $\ov{A}_2$ and
$f_1 (\ov{A}_2) \neq 0$
\\[1mm]
($ii$) $a=0$, $ f_0 (A_2)$ admits a strict local minimum at
the point $\ov{A}_2$ and
$f_1 (\ov{A}_2 + i  (\pi /2) \beta) \neq 0$.
\\[2mm]
Then, for all small $\e$ such that $\om_{\e} = (1/\sqrt{\e} ,
\beta \e^a)$ satisfies
$$
\om_{\e} \cdot {\bf k} \geq \frac{\gamma_\e}{|{\bf k}|^{\tau}},
\ \forall k \in {\bf Z}^n, k \neq 0
$$
for all $ I_0, I_0'$ with $\om_{\e} \cdot I_0= \om_{\e} \cdot I_0'$,
there is a heteroclinic orbit connecting the invariant tori
${\cal T}_{I_0}$ and
${\cal T}_{I_0'}$. In addition,
for all $ \eta > 0 $ small enough the ``diffusion time'' $T_d$ needed
to go from a $\eta$-neighbourhood of
${\cal T}_{I_0}$ to a $\eta$-neighbourhood of
${\cal T}_{I_0'}$ is
$O( | I_0 - I_0' | (\sqrt{\e}/\mu) e^{\pi/(2\sqrt{\e})} [(\gamma_\e)^{-1}
(\sqrt{\e} e^{\pi
/(2\sqrt{\e})})^{\tau} + |\ln \mu|]
+  | \ln ( \eta )|).$
\end{theorem}

\begin{pf}
It is enough to prove that, if ($i$) or ($ii$) is satisfied, then the
condition given in lemma \ref{lem:3ts} holds. The statement
is then a direct consequence of theorem \ref{thm:main}.

We first assume that condition ($i$) is satisfied. In what follows,
the notation $ u = O(v) $ means that $ | u | \leq C |v|$, where $ C $ is
a universal constant.
We have
\begin{eqnarray*} 
\Gamma_0 (\e, A_2)&= &\sum_{\kappa_2 \in {\bf Z}^{n-1}}
 \frac{2\pi \kappa_2 \cdot \beta \e^a }
{\sinh (\pi \kappa_2 \cdot \beta \e^a /2)} f_{0,\kappa_2}
e^{i\kappa_2 \cdot A_2} \\
&=& \sum_{\kappa_2 \in {\bf Z}^{n-1}} (4+O(\e^{2a} |\kappa_2|^2)) 
f_{0,\kappa_2} e^{i\kappa_2 \cdot A_2} \\
&=& 4 f_0 (A_2) + O \Big( \sum_{\kappa_2\in {\bf Z}^{n-1}} 
\e^{2a} |\kappa_2|^2 e^{-\kappa_2^+ \cdot \rho_2} \Big) \\
&=& 4f_0(A_2) + O ( \e^{2a} ).
\end{eqnarray*}
Moreover 
\begin{eqnarray*}
\Gamma_1 (\e, A_2)&= &\sum_{\kappa_2 \in {\bf Z}^{n-1}} \frac{2\pi (\e^{-1/2}
+  \kappa_2 \cdot \beta \e^a )}
{\sinh ((\pi/2)(\e^{-1/2 } +  \kappa_2 \cdot \beta \e^a ))} f_{1,\kappa_2}
e^{i\kappa_2 \cdot A_2} \\
&=&
\sum_{\kappa_2\in {\bf Z}^{n-1}} \frac{4\pi}{\sqrt{\e}}
e^{-(\pi/2)(\e^{-1/2} + \kappa_2 \cdot \beta \e^a )} 
\Big( 1+O(|\kappa_2| \e^{a+(1/2)}) \Big) f_{1,\kappa_2}
e^{i\kappa_2 \cdot A_2}
\\
&=& \frac{4\pi}{\sqrt{\e}}
e^{-(\pi/2)\e^{-1/2}} \Big[ f_1 (A_2) +  O \Big( 
\sum_{\kappa_2\in {\bf Z}^{n-1}}
e^{-\kappa_2^+ \cdot \rho_2 } 
( | e^{-(\pi/2)\kappa_2 \cdot \beta \e^a } -1|
+ |\kappa_2| \e^{a+(1/2)} e^{|\kappa_2 \cdot \beta| \e^a}) \Big) \Big]
\\
&=& \frac{4\pi}{\sqrt{\e}}
e^{-(\pi/2)\e^{-1/2}} \Big[ f_1 (A_2) +  O \Big( 
\sum_{\kappa_2\in {\bf Z}^{n-1}}
\exp(-\sum_{j=2}^n |k_j|(r_j-(\pi/2)|\beta_j|\e^a)) |\kappa_2| \e^a \Big) \Big]
\\
&=& \frac{4\pi}{\sqrt{\e}}
e^{-(\pi/2)\e^{-1/2}} \Big[ f_1 (A_2) +  O(  \e^a) \Big],
\end{eqnarray*}
provided $\e$ is small enough.
It is then clear that condition ($i$) implies that assumption of lemma
\ref{lem:3ts} holds.

We now assume that condition ($ii$) is satisfied.
As previously, we have
\begin{eqnarray*}
\Gamma_1 (\e, A_2)&= &
\sum_{\kappa_2\in {\bf Z}^{n-1}} \frac{4\pi}{\sqrt{\e}}
e^{-(\pi/2)(\e^{-1/2} + \kappa_2 \cdot \beta)} 
(1+O(|\kappa_2| \e^{1/2}) f_{1,\kappa_2}
e^{i\kappa_2 \cdot A_2}
\\
&=& \frac{4\pi}{\sqrt{\e}}
e^{-(\pi/2)\e^{-1/2}} \Big[ 
f_1 (A_2+i(\pi/2) \beta ) +  O \Big( \sum_{\kappa_2\in {\bf Z}^{n-1}}
\exp(-\sum_{j=2}^n |k_j|(r_j-(\pi/2)|\beta_j|) )   |\kappa_2| \e^{1/2} \Big)
\Big]
\\
&=& \frac{4\pi}{\sqrt{\e}}
e^{-(\pi/2)\e^{-1/2}} \Big[ f_1 (A_2+i(\pi/2) \b ) +  O( \sqrt{\e}) \Big].
\end{eqnarray*}

We observe also  that, if $ a=0 $, then
$ \Gamma_0 ( \e, A_2) $ is independent of $\e$.
It follows easily that condition ($ii$) implies that the assumption
of lemma \ref{lem:3ts} holds true.
\end{pf}

\begin{remark}
In many examples, condition ($ i $) or condition ($ii$)
is satisfied. However we need that
$ f_0 (A_2) $ and $f_1 (A_2) $ do not vanish everywhere,
see remark \ref{rem:tts}-(iv).
\end{remark}

\section{Appendix}

In the proof of the following lemmas we will  closely follow
the arguments developed in the papers \cite{BB}-\cite{BBS} to which
we refer for further details.
In the sequel the notation $u=O(v)$ (resp. $u=o(v)$) will
mean that there is a constant $C$ (resp. a function $\e(v)$) 
independent of anything except $f$ such that 
$|u|\leq C |v|$ (resp. $|u| \leq \e (v) |v|$ and $\lim_{v\to 0}
\e(v) =0$).
\\[1mm]
\begin{pfn}{\sc of lemma } \ref{lem:1bump}.
We first assume that $\theta=0$ and give the existence proof 
in $ [ 0 , + \infty ) $.
We are looking for a solution of (\ref{pendper}) in the form of
$ q = q_0 + w $ with $ w ( 0 ) = 0 $ and
$\lim_{t \to + \infty} w(t) = 0 $. The function $ w $
must  satisfy the equation
$$
- \ddot w + w = - \Big( \sin (q_0 + w) - \sin q_0
- w \Big) + \mu \sin ( q_0 + w ) f (\om t + A).
$$
Let
$$
{\bf X} = \Big\{ 
w( \cdot ) \in W^{1, \infty} ([0, +\infty )) \ \Big| \
||w||_1 := 
\sup_{t \in {\bf R}} 
\max (|w(t)|, |{\dot w}(t)|) \exp({\frac{| t |}{2}}) < + \infty \Big\}
$$
and 
$$
{\bf X}' =
\Big\{ 
w( \cdot ) \in L^{\infty} ([0, +\infty )) \ \Big| \ 
||w||_0 := \sup_{t \in {\bf R}} 
|w(t)| \exp({\frac{|t|}{2}}) < + \infty \Big\}.
$$
${\bf X}$ and ${\bf X}'$, endowed respectively with norms
$|| \ ||_1$ and $|| \ ||_0$, are Banach spaces.
Let ${\cal L}_0 $ be the linear operator which
assigns to $ h \in {\bf X}'$ the unique solution $ u = {\cal L}_0 h $ of the
problem:
$$
\left\{  \begin{array}{l}
- {\ddot u } + u = h \\
 u(0)= 0 \ ,  \
\lim_{t \to +\infty}  u(t) = 0.
\end{array}  \right.
$$
An explicit computation shows that, for $ t \in [0, +\infty ) $,
\be \label{eq:s}
 u ( t )  = ({\cal L}_0 h)(t) = 
\frac{1}{2}  \int_{0 }^{+\infty}
\Big( e^{-|t - s|}-e^{-(t + s)} \Big) h(s) ds.
\ee
As an easy consequence ${\cal L}_0$ sends ${\bf X}'$ into 
${\bf X}$ continuously.

We define the non-linear operator
$ H:{\bf R} \times {\bf R}^n  \times {\bf X} \to {\bf
X}$ by 
\be\label{eq:bans}
H ( \mu, A,  w ) :=
w - {\cal L}_0 \Big( - \Big( \sin ( q_0 + w) - \sin q_0
- w \Big) + \mu \sin ( q_0 + w ) f ( \om t + A) \Big).
\ee
$H$ is smooth, $2\pi {\bf Z}^n$-periodic w.r.t. $A$ and we have $H(0,A,0)=0$.
The unknown $ w $ must solve the equation $H(\mu, A, w ) = 0$.
We can apply the Implicit Function Theorem.
In fact, let us check that
$$
 \partial_w  H (0,A,0) :  \  W  \to  
W - {\cal L}_0 \Big[  (1-\cos q_0) W  \Big]
$$
is invertible. Since $\lim_{t\to \infty}(1-\cos q_0(t))=0$, 
$ \partial_w H (0,A,0) $
is  of the type``Identity + Compact'' and then it is sufficient to show that
it is injective. $W$ is in the kernel of
$ \partial_w H( 0, A,  0)  $ iff $W(0)=0$ and W satisfies in 
$( 0, + \infty)$ the equation
\be\label{eq:linea}
- \ddot{W} + \cos q_0 W = 0.
\ee
Multiplying by $ {\dot q}_0 $ in (\ref{eq:linea}) and
integrating over $[ 0, + \infty) $ by parts twice we obtain that  
$ {\dot W}(0 ) {\dot q}_0 (0) =0 $. 
Since $ {\dot q}_0 (0) \neq 0 $ we 
get also $ {\dot  W}(0 ) =0 $ and then $W = 0$.
Thus the kernel of $ \partial_w H( 0, A,  0)  $ is reduced to
$0$, and this operator is invertible.
We derive by the Implicit Function Theorem that there 
are $\rho_0>0$ and $\mu_0>0$ such that, for all 
$|\mu| < \mu_0$, for all $A\in {\bf R}^n$,  
the equation $H(\mu,A,w)=0$  has a unique
solution $w_A^{\mu}$ in ${\bf X}$ such that $||w^{\mu}_A||< \rho_0$.

Note that $\mu_0$ and $\rho_0$ may be chosen independent of 
$A$ (and of $\om$ too) because $\partial_w H(0,A,0)$ is independent
of $A$ and $\om$, $\partial_{\mu} H(0,A,0)$ is uniformly bounded,
and $\partial_w H(\mu,A,w)$ (resp. $\partial_{\mu} H(\mu,A,w)$)
tend to $\partial_w H(0,A,0)$ (resp. $\partial_{\mu} H(0,A,0)$)
as $(\mu,w) \to (0,0)$ uniformly in $(A,\om)$.

Since $H$ is smooth $w_A^{\mu}$ depends smoothly on $\mu$ 
and $A$ and $w_{A+2\pi k}^{\mu}=w_A^{\mu}$ by the $2\pi {\bf
Z}^n$-periodicity of $H$ w.r.t. $A$.
  By the properties of $\partial_{\mu} H$ mentioned
above, $||w_A^{\mu}||_1 = O(\mu)$.

In a similar way we can prove the existence and unicity of
${w'}_A^{\mu} : (-\infty,0] \to {\bf R}$ which satisfies analogous
properties over the interval $(-\infty,0]$.
We can define $q_{A,0}^{\mu}$ by 
$q_{A,0}^{\mu} (t)=q_0(t)+w_A^{\mu}(t)$ if $t\geq 0$,
$q_{A,0}^{\mu}=q_0(t)+{w'}_A^{\mu}(t)$ if $t<0$. This is
the unique function for which $(i)$, $(ii)$ (with $\theta=0$)
and $(iii)$ hold.

If $\theta \neq 0$, we observe that $q$ satisfies $(i)$ iff
$$
\left\{
\begin{array}{l}
- (T_{-\theta} q)'' + \sin (T_{-\theta} q) = \mu 
\sin(T_{-\theta} q) f(\om t+ A+\om \theta) \\
(T_{-\theta} q) (0)=\pi,
\end{array}
\right. 
$$
where $T_{-\theta} q (t)=q(t+\theta)$. Hence there is a unique 
$q_{A,\theta}^{\mu}$ which satisfies $(i)$, $(ii)$, defined
by $q_{A,\theta}^{\mu}=T_{\theta} q_{A+\om \teta,0}^{\mu}$, i.e.
$q_{A,\theta}^{\mu}(t)=q_{A+\om \theta,0}^{\mu} (t-\theta)$;
$(iii)$ and $(iv)$ clearly hold. The regularity of
$q_{A,\theta}^{\mu}$ w.r.t. $A,\theta,\mu$ is a consequence of
the regularity of $w_A^{\mu}$ and ${w'}_A^{\mu}$ w.r.t. $A$ and
$\mu$.
$(v)$ follows from
$$
\partial_A w_A^{\mu} = - \Big [ \partial_{w} H
(\mu,A,w_A^{\mu}) \Big ]^{-1}  \partial_A H(\mu,A,w^A_{\mu})
$$
provided we can justify that
$ ||\partial_A H(\mu,A,w_A^{\mu})||_1 = O(\mu)$, $\quad
||\om \cdot \partial_A H(\mu,A,w_A^{\mu})||_1 =O(\mu).$
The second bound (uniform in $\om$) is not so obvious.
We just point out that
$$
\om \cdot \partial_A H(\mu,A,w_A^{\mu}) = -{\cal L}_0 (\mu \sin(q_0 +
w_A^{\mu}) \frac{d}{dt} f(\om t +A)
$$
and that we can use the ``regularizing'' properties of ${\cal L}_0$.
\end{pfn}

\begin{pfn}{\sc of lemma} \ref{lem:kheter}.
We give the proof in the interval $[\teta_1 , \teta_2 ]$.
We may assume without loss of generality that $\theta_1=0$
since, by the remark at the end of the proof of lemma
\ref{lem:1bump}, a translation of the time by  $-\teta_1$
amounts to adding $\om \theta_1$ to $A$. For simplicity
of notations, we shall write $\theta_2=\theta$.

We are  looking for a solution $ q = q^*_{0, \teta } + w $ of
(\ref{pendper}) over $(0,\theta)$
with $ w ( 0 ) = w ( \teta ) = 0 $, where
$q^*_{ 0, \teta } $ is the following
smooth ``approximate solution''
$$
\begin{array}{rcl}
q^*_{ 0, \teta } (t)= 
\end{array}
\left\{
\begin{array}{rcl}
q^\mu_{A,0} (t) \; {\rm if} \; \
t \in  ( 0 ,   \teta / 2 - 1 ), \\
r^*_\teta (t) \; {\rm if} \; \ t \in
 [\teta / 2 - 1,   \teta / 2 + 1]\\
2 \pi + q^\mu_{ A, \teta } (t) \; {\rm if} \; t \in
(\teta / 2 + 1, \teta),
\end{array}\right.
$$
where
$$
 r^*_{\teta} (t)=(1-R(t-\theta/2)) q^\mu_{A,0} (t)  +
R(t-\theta/2) (q^\mu_{A,\theta} (t) + 2\pi),
$$
and 
$R : {\bf R} \to [0,1]$ is a $C^{\infty}$ function such that 
$R(s)=0$ if $s\leq -1$, $R(s)=1$ if $s\geq 1$.
Let $ {\cal L}_{0, \teta} $ be the linear operator which
assigns to $ h \in L^{\infty}([0,\theta]) $ the unique solution
$ u = {\cal L}_{0, \teta} h $ of the problem:
\be\label{eq:twoz}
\left\{
\begin{array}{l}
- {\ddot u } + u = h \\
 u( 0)= 0 \ , \
u( \teta ) = 0.
\end{array}
\right.
\ee
An explicit computation shows that for $t \in [0, \teta ] $
the solution $ u $ of (\ref{eq:twoz}) is given by
$$
u (t)  =   \frac{1}{\sinh(\teta)} \Big[
\int_0^{t}
h (s ) \sinh ( s ) \sinh (\teta - t) \ ds  \
 +  \
\int_{t}^{\teta} h (s) \sinh (\teta-s) \sinh (t) \ ds \Big].
$$
Note that $ {\cal L}_{0, \teta}$  sends $L^{\infty} ([0,\theta])$
into $W^{1,\infty} ([0,\theta])$   ($W^{2,\infty} ([0,\theta])$ 
in fact) and that there is a constant $K$
independent of $\theta$ such that $|| {\cal L}_{0, \teta}
W||_{1,\infty} \leq K ||W||_{\infty}$, where
$||\ ||_{\infty}$ denotes the infty norm in $[0,\theta]$ and
$||W||_{1,\infty}:=||W||_{\infty}+||\dot{W}||_{\infty}$.

We define the smooth non-linear operator   
$ H^{\theta} : {\bf R} \times {\bf R}^n  \times 
W^{1, \infty } ([0, \teta]) \to W^{1, \infty } ([ 0, \teta ])$ by
$$
H^{\theta}( \mu, A,  w) :=  
w - {\cal L}_{0, \teta} 
\Big( - \Big( \sin (q^*_{0, \teta} + w)
- {\ddot q}^*_{0, \teta} - w \Big) + \mu 
\sin ( q^*_{0, \teta} + w ) f (\om t + A) \Big).
$$
We immediately remark for further purpose that
\be \label{D2H}
||\partial^2_{ww} H^{\theta} (\mu,A,w) [W, W]|| 
= O(||W||^2_{\infty}).
\ee
Moreover, by lemma \ref{lem:1bump}-$(i)$ and 
 the definition of $ q^*_{0,\theta} $,
$|| - \sin q^*_{0, \teta} + {\ddot q}^*_{0, \teta}  + \mu 
\sin ( q^*_{0, \teta}  ) f (\om t + A) \Big)  ||_{\infty} =
O(\exp(-\theta/2))$ hence
\be \label{upHH}
||H^{\theta}(\mu,A,0)||_{1,\infty} =O(\exp(-\theta/2)).
\ee

$q^*_{0, \teta} + w$ is a solution of (\ref{pendper})
with the appropriate boundary conditions iff
$H^{\theta} (\mu, A, w) = 0  $
 
We shall show that there
exist $ \ov{C}, \ov{L}, \ov{\mu} > 0 $ such that
$\forall \theta > \ov{L} $,  for all
$ | \mu | < \ov{\mu} $, for all $A$ and $\om$,
$\partial_w H^{\theta}(\mu, A, 0)$ is invertible and
\be
\Big| \Big| \left(\partial_w H^{\theta}
(\mu, A, 0)\right)^{-1} \Big| \Big| \leq \ov{C}.
\ee
Since $\partial_w H^{\theta} (\mu, A, 0)$ is of the type 
``Id + Compact'', it is enough to prove that 
$$
\forall W \in W^{1,\infty} ([0,\theta]) \quad
|| \partial_w H^{\theta} (\mu, A, 0)  W||_{1,\infty}
\geq \frac{1}{\ov{C}} ||W||_{1,\infty}.
$$
We shall just sketch the proof of this assertion (see
also lemma 2 of \cite{BB}). Arguing by contradiction, we assume
that there are sequences $(\mu_n) \to 0$, $(\theta_n)\to\infty$,
 $(A_n)$, $(\om_n)$, $(W_n)$ such that
$W_n \in W^{1,\infty} ([0,\theta_n]), ||W_n||_{1,\infty}=1$,
\be \label{H,w}
||\partial_w H^{\theta_n} (\mu_n, A_n, 0) W_n||_{1,\infty} \to 0.
\ee

Let $\xi_n \in [0,\theta_n]$ be such that
$
m_n:={\rm max}_{t\in [0,\theta_n]} |W_n(t)|=W_n(\xi_n).
$  
By (\ref{H,w}) and the properties of ${\cal L}_{0,\theta}$,
$||W_n||_{1,\infty}=O(m_n)$. Hence 
$\liminf(m_n) >0$. Taking a subsequence, 
we may assume that  $(\xi_n)$ is  bounded  or 
$(\theta_n -\xi_n)$ is bounded or ($(\xi_n)\to \infty$ and 
$(\theta_n -\xi_n) \to \infty$).

In the first case, still up to a subsequence 
 $W_n \to W \neq 0$ uniformly in compact subsets of
$[0,\infty)$. Taking limits in (\ref{H,w}) we obtain that
$W(0)=0$, $-\ddot{W} + \cos q_0  W=0$, which 
contradicts $W\neq 0$. The second case can be dealt with
similarly. In the third case, up to a subsequence,
$W_n ( \cdot + \xi_n) \to W \neq 0$ uniformly in compact subsets
of ${\bf R}$, with $|W(t)| \leq |W(0)|$ for all $t\in{\bf R}$.
Taking limits in (\ref{H,w}), we obtain that
$-\ddot{W}+W=0 $ over ${\bf R}$, which contradicts
$W\neq 0$ bounded.  

From now we shall assume that $| \mu | < \mu_0 \leq \ov{\mu}$,
$\theta > \ov{L}$.   Let
$$ R^{\theta}(\mu, A,  w) = H^{\theta}(\mu, A , w) 
- H^{\theta}(\mu, A,  0) - 
\partial_w H^{\theta}(\mu, A, 0) w . $$
By the previous assertion, 
$$
H^{\theta} ( \mu, A,  w) = 0 \quad \Leftrightarrow \quad
w = - \left(\partial_w H^{\theta}
(\mu, A, , 0) \right)^{-1}
H^{\theta} ( \mu, A, 0)  
- \left(\partial_w H^{\theta}(\mu, A,  0)
\right)^{-1} R^{\theta}(\mu, A, w )
:= F^{\theta}_{\mu, A}( w ).
$$
We just have to show that
$ F^{\theta}_{\mu, A} $ is a
contraction in some ball $ B(0,\rho)  \subset W^{1,\infty} ([0,\theta])$.
For this, we derive from (\ref{D2H}) and (\ref{upHH}) in a
standard way that,
for all $ || w ||_{1,\infty}, || w' ||_{1,\infty} \leq \rho $,
$|\mu| < \ov{\mu}$, $\theta>\ov{L}$ there holds
\begin{equation} \label{eq:disu}
\| F^{\theta}_{\mu, A} ( w ) \|_{1,\infty} = O(
\exp({-  \theta /2 }) +   \rho^2) \quad  ;  \quad
\| F^{\theta}_{\mu, A} ( w ) -
   F^{\theta}_{\mu, A} ( w') \| =O( \rho || w' - w ||).
\end{equation}

We can deduce that $F^{\theta}_{\mu,A}$ is a contraction in
$\ov{B}(0,\rho)$, with
$ \rho =  C \exp({-  \theta/2})$,  for some constant $C$,
provided that $\theta > \ov{L}$,  $\ov{L}$ large enough.
Applying the Contraction Mapping Theorem we conclude that there is
a unique solution
$|| w^L_\mu (A,\teta ) ||_{1,\infty} \leq  C \exp({ -  \theta /2 }) $
of  the equation $H^{\theta}_{\mu,A} (w)=0$.
Note that by (\ref{eq:disu}) uniqueness holds in $B(0,\rho_0)$
for some $\rho_0>0$ independent of $\theta$. The regularity
of the solutions in $(A, \teta, \mu)$
follows like in \cite{BB}. 
\end{pfn}

\begin{pfn}{\sc of lemma} \ref{lem:prdi}.
Let us  consider the function
$ H: {\bf R} \times {\bf T}^n \times {\bf R} \times {\bf R} \to {\bf R}$
defined by
$$
H ( \mu, A, \teta, l ) = Q_{A, \teta}^\mu ( \teta + l ) - \pi
$$
The unknown $ l_\mu (A, \teta) $ can be implicitely defined
by the equation $ H ( \mu, A, \teta, l ) = 0 $.
We have  $ H ( 0, A, \teta, 0 ) = 0 $ and
$$
\partial_l  H ( 0, A, \teta, 0 ) =  {\dot q}_0 ( 0 ) \neq 0.
$$
Hence by the Implicit function theorem,
for $\mu$ small enough (independently of $A,\theta,\om$
because $\partial_l H$ and $\partial_{\mu} H$ are
continuous uniformly in $A,\theta,\om$), there exists a unique smooth solution
$ l_\mu (A, \teta) = O ( \mu )$ of $ H ( \mu, A, \teta, l ) = 0 $.
Moreover, by the uniform estimates in $A$ and $\om$ that 
we can obtain for $\partial_A Q^{\mu}_{A,\theta}$, 
$\om \cdot\partial_A Q^{\mu}_{A,\theta}$, there holds
$|\nabla l_{\mu} (A, \theta)|=O(\mu)$.
\end{pfn}

\begin{pfn}{\sc of lemma} \ref{lem:apfuf}.
The first step  is to prove that
\be\label{eq:primaes}
\max \Big(|q^\mu_{A,\teta + l_\mu (A, \teta)}(t) - Q_{A,\teta}^\mu (t)|,
|{\dot q}^\mu_{A, \teta + l_\mu (A, \teta)}(t) - {\dot Q}_{A,\teta}^\mu (t)|
\Big) \leq K_0 | \partial_\teta \wtilde{F}_\mu (A, \teta) |
\exp({- \frac{ |t - \teta |}{2}}), \forall t \in {\bf R}.
\ee
We have
$q^{\mu}_{A,\theta+l_{\mu}(A,\theta)}=T_{\theta+l_{\mu}(A,\theta)}
q^{\mu}_{A',0}$; $\ Q^{\mu}_{A,\theta}=T_{\theta+l_{\mu}(A,\theta)}
Q^{\mu}_{A', \theta'}$, where
$A'=A+\om(\theta+l_{\mu}(A,\theta))$ and
$\theta' = - l_{\mu} (A,\theta)$. So it
is enough to prove the estimate for
$ \
w:=Q^{\mu}_{A', \theta'}-q^{\mu}_{A',0}.$

Note that $Q^{\mu}_{A', \theta'} (0)=q^{\mu}_{A', 0}(0)=\pi$. 
So $Q^{\mu}_{A', \theta'}-q_0$, $q^{\mu}_{A', 0}-q_0$ belong to
${\bf X}$  and satisfy
$$
H(\mu,A',q^{\mu}_{A', 0}-q_0)=0 \quad ; \quad 
H(\mu, A', Q^{\mu}_{A', \theta'}-q_0)= \alpha^{\mu}_{A',\theta'}
{\cal L}_0 (\psi_{\theta'}),
$$
where ${\bf X}, H$ and $ {\cal L}_0$ are defined in the proof of lemma
\ref{lem:1bump}. Therefore
$$
\alpha^{\mu}_{A',\theta'} {\cal L}_0 (\psi_{\theta'})=
\partial_w H(\mu, A', q^{\mu}_{A', 0}-q_0) (Q^{\mu}_{A', \theta'}-
q^{\mu}_{A', 0}) + o(||Q^{\mu}_{A', \theta'}-
q^{\mu}_{A', 0}||_1).
$$
Moreover $||Q^{\mu}_{A', \theta'}-q_0||_1 + || q^{\mu}_{A', 0}
- q_0||_1 =O(\mu)$. Hence, by the properties of $H$
mentioned in the proof of lemma \ref{lem:1bump}
(in particular the fact that $\partial_w H(0, A', 0)$ is
invertible) we obtain, for $\mu$ small enough the following
bound :
$$
||Q^{\mu}_{A', \theta'} - q^{\mu}_{A', 0}||_1 = O(|\a^{\mu}_{A',\theta'}|).
$$
Since
$| \a^\mu_{A', \teta'} | = O(|\partial_\teta \wtilde{F}_\mu (A, \teta)|)$
we deduce  estimate (\ref{eq:primaes}).


We can now estimate $\wtilde{F}_\mu ( A, \teta) - V_\mu (A,
\teta)$. For $q\in {\bf X}$ let
$$
{\cal G}^{\mu}_{A'} (q) =\int_{\bf R}
{\cal L}_{\mu,A'} (q,\dot{q}, t) \ dt
=\int_{\bf R} \frac{1}{2} (\dot{q}(t))^2
+ (1-\cos(q(t))) + \mu (\cos(q(t))-1) f(A'+\om t) \ dt.
$$
By standard arguments, ${\cal G}^{\mu}_{A'} : {\bf X} \to {\bf R}$
is smooth. Moreover for $ |\mu | <\mu_0$,
$$
D^2 {\cal G}^{\mu}_{A'} (q) [ w , w]= \int_{\bf R}
 \dot{w}^2 + \cos (q) w^2  - \mu \cos (q) w^2 
f(A'+\om t) \ dt = O(||w||^2_1).
$$
By the definition of 
$q^{\mu}_{A', 0}$ ($(i)$ in lemma \ref{lem:1bump}), we easily 
obtain with an integration by parts that
$D{\cal G}^{\mu}_{A' }(q^{\mu}_{A', 0})w=0$ for all $w\in {\bf X} $ such that
$w(0)=0$. Therefore 
$$
{\cal G}^{\mu}_{A'}(q^{\mu}_{A', 0}+w)=
{\cal G}^{\mu}_{A'}(q^{\mu}_{A', 0}) + O(||w||^2_1)
$$
for all $w\in {\bf X}$ such that $w(0)=0$. Hence
since $(Q^{\mu}_{A', \theta'}-q^{\mu}_{A', 0})(0)=0$, 
\begin{eqnarray*}
\wtilde{F}_\mu ( A, \teta)-V_{\mu}(A,\theta)   &=&
 \wtilde{F}_\mu ( A, \teta + l_{\mu} (A,\theta) + \teta')-
F_\mu ( A, \teta + l_{\mu}(A,\teta))   
   = \wtilde{F}_\mu ( A', \teta') - F_\mu ( A', 0)    \\
   &=& {\cal G}^{\mu}_{A'} (Q^{\mu}_{A', \theta'})-
{\cal G}^{\mu}_{A'} (q^{\mu}_{A', 0})
 = O(||Q^{\mu}_{A', \theta'} - q^{\mu}_{A', 0}||^2_1).               
\end{eqnarray*}
We obtain by (\ref{eq:primaes}) that
$$
|\wtilde{F}_\mu ( A, \teta)-V_{\mu}(A,\theta) | 
\leq   C_4 \Big( \partial_\teta \wtilde{F}_\mu (A, \teta) \Big)^2
$$
for some positive constant $ C_4 $.
\end{pfn}

\noindent
{\it Massimiliano Berti, S.I.S.S.A., Via Beirut 2-4,
34014, Trieste, Italy, berti@sissa.it}.
\\[2mm]
{\it Philippe Bolle,
D\'epartement de math\'ematiques, Universit\'e
d'Avignon, 33, rue Louis Pasteur, 84000 Avignon, France,
philippe.bolle@univ-avignon.fr}

\end{document}

%% file: macro.tex
\newenvironment{pf}{\noindent{\sc Proof}.\enspace}{\rule{2mm}{2mm}\medskip}
\newenvironment{pfn}{\noindent{\sc Proof} \enspace}{\rule{2mm}{2mm}\medskip}
\newtheorem{theorem}{Theorem}[section]
\newtheorem{lemma}{Lemma}[section]

\newtheorem{condition}{Condition}[section]
\newtheorem{remark}{Remark}[section]
\newtheorem{remarks}{Remark}[section]
\newtheorem{definition}{Definition}[section]
\newcommand{\be}{\begin{equation}}
\newcommand{\ee}{\end{equation}}

\newcommand{\teta}{\theta}

\newcommand{\om}{\omega}

\newcommand{\ov}{\overline}

\newcommand{\wtilde}{\widetilde}

\renewcommand{\a }{\alpha }
\renewcommand{\b }{\beta }
\newcommand{\s }{\sigma }
\renewcommand{\d }{\delta }

\newcommand{\e }{\varepsilon }
\newcommand{\g }{\gamma}

\newcommand{\vphi}{\varphi }

\renewcommand{\t }{\tau }

